\documentclass[12pt,reqno]{amsart}
\usepackage{amsfonts}
\usepackage{amscd}
\usepackage{amssymb}
\usepackage{amsfonts}
\usepackage{amscd}
\usepackage{amsmath} 
\usepackage{mathrsfs}
\usepackage{enumerate}
\usepackage{float} 
\usepackage[pdftex]{graphicx}
\allowdisplaybreaks

\usepackage{hyperref}

\date{\today}

\newtheorem{thm}{Theorem}[section]

\newtheorem{lem}[thm]{Lemma}

\newtheorem{prop}[thm]{Proposition}
\theoremstyle{definition}

\theoremstyle{remark}
\newtheorem{rem}[thm]{\bf{Remark}}

\numberwithin{equation}{section}

\newcommand{\la}{\lambda}

\newcommand{\gf}{\mathfrak{g}}
\newcommand{\kf}{\mathfrak{k}}
\newcommand{\pf}{\mathfrak{p}}
\newcommand{\af}{\mathfrak{a}}
\newcommand{\ad}{\text{ad}}
\newcommand{\Ad}{\text{Ad}}

\renewcommand{\Re}{\operatorname{Re}}
\renewcommand{\Im}{\operatorname{Im}}
 
\usepackage{color}

\title[Rellich-type estimates]
{On Rellich-type asymptotics for eigenfunctions on rank one symmetric spaces of  noncompact type }

\author[Ganguly]{ Pritam Ganguly}

\address[Pritam Ganguly]{Stat-Math Unit,
	Indian Statistical Institute, Kolkata,
	BT Road, Baranagar, Kolkata 700108}

\email{pritam1995.pg@gmail.com}

\date{}

\begin{document}
	
	\maketitle
	
	\begin{abstract}
		We study eigenfunctions of the Laplace--Beltrami operator \(\Delta_X\) in exterior domains \(\Omega\) of rank-one Riemannian symmetric spaces of noncompact type \(X\), a class that includes all hyperbolic spaces. Extending the classical \(L^2\) Rellich theorem for the Euclidean Laplacian, we analyze the asymptotic behaviour and \(L^p\)-integrability of solutions to the Helmholtz equation
		\[
		\Delta_X f + (\lambda^2 + \rho^2) f = 0 \quad \text{in } \Omega,
		\]
		where \(\lambda \in \mathbb{C}\setminus i\mathbb{Z}\) and \(\rho\) denotes the half-sum of positive roots.
		
		We establish sharp Rellich-type quantitative \(L^p\)-growth estimates in geodesic annuli, which yield the nonexistence of nontrivial \(L^p(\Omega)\)-solutions in the optimal range \(1 \leq p \leq 2\) for spectral parameters satisfying \(|\Im(\lambda)| \leq (2/p - 1)\rho\). For non-real spectral parameters, we further obtain refined Rellich-type uniqueness results under weak \(L^p\)-assumptions. As a by-product, we also prove a Rellich-type uniqueness theorem in terms of Hardy-type norms.
		
		Our results provide a geometric extension of the Euclidean Rellich theorem, highlighting the role of exponential volume growth and the \(p\)-dependence of the \(L^p\)-spectrum of \(\Delta_X\) in producing genuinely non-Euclidean spectral phenomena.
	\end{abstract}
	
	\section{Introduction and main results}
	The study of the asymptotic behavior of eigenfunctions of differential operators, such as the Laplacian and Schr\"odinger operators, is a central theme in spectral theory and mathematical physics. Such results are particularly valuable for ruling out the presence of eigenvalues embedded in the continuous spectrum of the time-independent Schr\"odinger operator, a fact of immense importance in quantum mechanics.  A classical result in this direction, describing an asymptotic estimate of eigenfunctions at infinity, is the following theorem, originally established by Rellich \cite[Satz 1]{Rel} in 1943.
	\begin{thm}[Rellich]
		\label{rel-org}
		Let $\lambda>0$, and $\Omega:=\{x\in \mathbb{R}^n: |x|>R_0\}$ for some $R_0>0.$ Assume that $0\neq f\in C^2(\Omega)$ satisfies the Helmholtz equation 
		\begin{equation}
			\label{helmholtz-eu}
			\Delta_{\mathbb{R}^n}f+\lambda f=0 \quad\quad \text{in} ~~\Omega.
		\end{equation}
		Then there exists $R_1>R_0$ sufficiently large, and a constant $C$ (depending on $\lambda, f$, and $n$) such that for all $R>R_1$ one has 
		\begin{equation}
			\label{rellich-asymp-euclidean}
			\int_{R<|x|<2R} |f(x)|^2~dx\geq C R.
		\end{equation}
	\end{thm}
	As is well known (see \cite[Anwendungen~2, p.~59]{Rel}), a direct consequence of \eqref{rellich-asymp-euclidean} is that the Helmholtz equation~\eqref{helmholtz-eu} 
	admits no nontrivial $L^2$–solutions in $\Omega$. 
	We must point out that establishing such a result on exterior domains (complements of bounded regions) 
	is considerably more delicate, whereas in the full space $\mathbb{R}^n$ the nonexistence of 
	$L^2$–solutions to $\Delta_{\mathbb{R}^n} f + \lambda f = 0$ with $\lambda>0$ follows immediately from the fact that 
	the Fourier transform provides a unitary equivalence between the Laplacian and 
	multiplication by $|\xi|^2$ on $L^2(\mathbb{R}^n)$.
	Moreover, Rellich’s $L^2$-asymptotic estimate, combined with unique continuation of Jerison and Kenig \cite{JK}, implies the absence of positive eigenvalues for the Schrödinger operator $H=-\Delta_{\mathbb{R}^n}+V(x)$ when the potential $V$ is compactly supported. The importance of the compactly supported case becomes clearer when one recalls the work of von Neumann and Wigner, who showed that highly oscillatory potentials can produce positive eigenvalues for the one-dimensional operator $H=-\tfrac{d^2}{dx^2}+V(x)$ (See \cite[p.223]{RS}). Later, in 1959, Kato~\cite{Kato} extended this phenomenon to all dimensions $n \ge 2$, 
	constructing potentials $V$ with $(1+|x|)V(x)\in L^\infty(\mathbb{R}^n)$ 
	for which the associated Schr\"odinger operator admits positive eigenvalues. 
	The subtle nature and usefulness of Rellich’s theorem subsequently inspired extensive research. 
	Littman~\cite{L} generalized Rellich’s result to partial differential equations with constant coefficients, 
	while further developments for Schr\"odinger operators were due to Kato~\cite{Kato}. 
	Building on this, Agmon~\cite{A} and Simon~\cite{S} extended Kato’s framework to potentials of the form 
	$V = V_1 + V_2$, where $V_1$ satisfies Kato’s condition and $V_2$ is a real-valued long-range term. 
	In 2003, Ionescu and Jerison~\cite{IJ} proved the absence of positive eigenvalues for $H$ 
	under suitable $L^q$-decay assumptions on $V$, removing the earlier compact-support restriction. 
	More recently, Rellich-type theorems have been established for generalized oscillators on $\mathbb{R}^n$, 
	characterizing the growth of eigenfunctions~\cite{T}, and for discrete Schr\"odinger operators~\cite{IM}.

	Returning to Rellich’s paper \cite{Rel}, he also noted another striking consequence of \eqref{rellich-asymp-euclidean}: if a solution to \eqref{helmholtz-eu} satisfies the size estimate 
	\begin{equation} \label{rellich-size-est-eu}  
		\lim_{|x|\to\infty} |x|^{\frac{n-1}{2}} |f(x)| = 0,  
	\end{equation}  
	then necessarily $f=0$. In other words, nontrivial eigenfunctions with positive eigenvalue cannot decay arbitrarily fast at infinity, unlike harmonic functions ($\lambda=0$), which may vanish arbitrarily fast even at infinity.

	Banerjee and Garofalo~\cite{BG1} recently revisited Rellich-type results, beginning with the observation that the size estimate~\eqref{rellich-size-est-eu} is, in fact, sharp.  
	Indeed, recall that the standard {Euclidean spherical functions}, given by  
	\begin{equation} 
		\label{sph-eu}  
		\phi_{\lambda}(x) 
		:= \int_{S^{n-1}} e^{-i \sqrt{\lambda}\, x \cdot \omega} \, \mathrm{d}\sigma(\omega)
		= c_{n} \frac{J_{n/2 - 1}(\lambda |x|)}{(\lambda |x|)^{n/2 - 1}}, 
		\qquad x \in \mathbb{R}^n,
	\end{equation}
	are well known to be radial solutions of the equation  
	$
	\Delta_{\mathbb{R}^n} f + \lambda f = 0
	$ for $\lambda>0.$
	Using the following well-known asymptotic behavior of the Bessel function \( J_{\nu} \) for \( \nu > -\tfrac{1}{2} \) (see, e.g., \cite[B.8]{Grafakos}),  
	\begin{equation} \label{eq:bessel-asymp}
		J_{\nu}(r) 
		= r^{-1/2} \!\left( \frac{2}{\pi} \right)^{1/2} 
		\cos\!\left( r - \frac{\pi \nu}{2} - \frac{\pi}{4} \right)
		+ R_{\nu}(r), 
		\qquad |R_{\nu}(r)| \leq c_{\nu} r^{-3/2},
		\quad ( r > 1),
	\end{equation}
	one observes that \( \phi_{\lambda} \) fails to satisfy~\eqref{rellich-size-est-eu}.  
	Furthermore, the same asymptotic behavior \eqref{eq:bessel-asymp} implies that  
	$
	\phi_{\lambda} \in L^{p}(\Omega) 
	\quad \text{if and only if} \quad 
	p > \frac{2n}{n - 1}.
	$

	Motivated by this observation, Banerjee and Garofalo investigated whether Rellich’s theorem admits an $L^p$ extension for 
	$p \leq \tfrac{2n}{n-1}$; equivalently, whether nontrivial solutions of~\eqref{helmholtz-eu} can exist in $L^p(\Omega)$. 
	Their main result, proved in~\cite{BG1}, is as follows.  
	
	\begin{thm}[Banerjee--Garofalo \cite{BG1}]
		\label{BG1-thm}
		Let $\lambda$ and $\Omega$ be as in Theorem~\ref{rel-org}.  
		Assume that $f$ is a solution of~\eqref{helmholtz-eu}. 
		If $f\in L^p(\Omega^*)$ with $\Omega^*:=\{x\in \mathbb{R}^n : |x|>R_0+4\}$ for some 
		$0 < p \leq \tfrac{2n}{n-1}$, then $f \equiv 0$ in $\Omega$.
	\end{thm}  
	
	As shown in~\cite{BG1}, this $L^p$ version follows from Rellich’s $L^2$ asymptotic estimate~\eqref{rellich-asymp-euclidean}.  In fact,
	for $p \in \left(2, \tfrac{2n}{n-1}\right]$, use of  H\"older’s inequality yields an $L^p$ analogue of~\eqref{rellich-asymp-euclidean}, 
	which immediately implies that no nontrivial $L^p(\Omega^*)$ solutions can exist. 
	For $0 < p < 2$, the proof instead relies on a Moser-type submean value estimate. 
	
	The $L^{p}$ extension of Rellich’s theorem due to Banerjee and Garofalo represents an 
	important development in the study of Rellich-type problems beyond the classical 
	$L^{2}$ framework. 
	As a continuation of this line of work, Banerjee and Garofalo~\cite{BG2} established 
	Rellich-type inequalities for the Baouendi–Grushin operators, using Carleman estimates 
	and weighted integral inequalities adapted to their subelliptic structure. 
	In a subsequent paper~\cite{BG3}, they extended Theorem~\ref{BG1-thm} to uniformly 
	elliptic divergence-form operators with asymptotically flat coefficients, covering 
	the range $0<p<\frac{2n}{n-1}$ while leaving the endpoint case open. 
	These contributions fit naturally into a broader program aimed at understanding 
	Rellich-type phenomena for a wider class of differential operators.

	In the present work, we investigate analogous questions in the setting of rank-one 
	Riemannian symmetric spaces of noncompact type, a class that includes the real and 
	complex hyperbolic spaces as fundamental examples, where exponential volume growth 
	and curvature play a significant role in determining the behaviour of 
	eigenfunctions.
	We commence by introducing the necessary notational framework to capture the intricacies of this setting and proceed to state the main results of the paper. Any notation or concept not explicitly explained at this stage will be clarified in Section~\ref{sec:rss}.

	Let \(X = G/K\) denote a rank-one Riemannian symmetric space of  noncompact type. Here, \(G\) is a connected,  noncompact, semisimple Lie group of real rank one with finite center, and \(K\) is a maximal compact subgroup of \(G\). We denote by \(d\) the geodesic distance function on \(X\), and by \(\Delta_{X}\) the Laplace–Beltrami operator associated with the \(G\)-invariant Riemannian metric on $X$. It is well known that, for any \(x \in X\), the volume of a geodesic ball \(B(x,r)\) of radius \(r\) satisfies
	\[
	|B(x,r)| \asymp e^{2\rho r}, \quad\quad (r > 1)
	\]
	where \(\rho > 0\) is a constant determined by the structure of \(X\) (see \eqref{rho-def}). This exponential volume growth plays a decisive role in the analysis on \(X\). 
	
	For $\lambda\in \mathbb{C}$, the elementary spherical functions given by 
	\begin{equation}
		\label{def:sphfn}
		\varphi_{\lambda}(x)=\int_K e^{-(i\lambda+\rho)H(x^{-1}k)}dk\quad\quad(x\in X)
	\end{equation}
	are important eigenfunctions of $\Delta_{X}$ which satisfies $\Delta_{X}\varphi_{\lambda}+(\la^2+\rho^2)\varphi_{\lambda}=0.$ For $\la\in\mathbb{R}$, the objects $e^{-(i\lambda+\rho)H(x^{-1}k)}$ are analogous to the functions $x\mapsto e^{-i\lambda x\cdot\omega}$ on $\mathbb{R}^n$, and the above representation \eqref{def:sphfn} can be compared with \eqref{sph-eu}, the situation in the Euclidean context.  This, in view of the self-adjointness, as consequence, implies that the $L^2$-spectrum  of $-\Delta_X$ is given by $S_2(-\Delta_{X})=[\rho^2, \infty).$ However, it is well known that for \(1 \leq p < \infty\), the \(L^{p}\)-spectrum of the Laplace--Beltrami operator \(-\Delta_{X}\) is given by a parabolic neighbourhood of the half-line \([\rho^{2}, \infty)\), namely (see \cite{LR}, \cite[Proposition~2.2]{Taylor})
	\begin{equation}
		\label{lp-spec}
		S_p(-\Delta_{X}) = \{ z^{2} + \rho^{2} : |\Im z| \leq \bigl| \tfrac{2}{p} - 1 \bigr| \rho \}.
	\end{equation}
	This description reveals a striking feature: in contrast with the Euclidean Laplacian, the \(L^{p}\)-spectrum of \(-\Delta_{X}\) depends in an essential way on \(p\). Such dependence reflects a deeper geometric dichotomy between spaces of polynomial volume growth and those of exponential volume growth (see \cite{Hu2, Hu3}).
	
	In the Euclidean case \(\mathbb{R}^{n}\), Rellich-type asymptotics in \(L^{p}\) can be deduced from the corresponding \(L^{2}\) theory (cf. \cite{BG1}). On symmetric spaces, however, this mechanism breaks down, and one is naturally led to expect genuinely new \(L^{p}\)-phenomena. This is particularly evident when considering eigenvalues with nonzero imaginary part, which have no analogue in the Euclidean setting.
	
	Furthermore, by the well-known estimates for spherical functions (see \eqref{phi-la-est}), one has \(\varphi_{\lambda} \in L^{p}(X)\) whenever \(p > 2\) and \(|\Im \lambda| < (1 - 2/p)\rho\) (See \eqref{phi-la-lp-est}). In particular, for \(p > 2\), the interior of \(S_p(-\Delta_{X})\) consists entirely of eigenvalues. 
	
	Motivated by these considerations, we investigate the asymptotic behaviour of eigenfunctions of the Laplace--Beltrami operator outside a compact obstacle. Our main result establishes a Rellich-type asymptotic in this setting, revealing new features that are intrinsic to the \(L^{p}\) framework on symmetric spaces.
	\begin{thm}
		\label{rellich-asymptotic-rss}
		Let $\lambda \in \mathbb{C}\setminus i\mathbb{Z}$, and set
		$
		\Omega := \{x\in X : d(o,x)>R_{0}\},
		$
		for some $R_0>0.$
		Suppose $0\neq f \in C^{2}(\Omega)$ satisfies the Helmholtz equation
		\begin{equation}
			\label{helmrss-intro}
			\Delta_{X}f + (\lambda^{2}+\rho^{2})f = 0 
			\qquad \text{in }\Omega.
		\end{equation}
		Then there exists $R_{1}>R_{0}$ (sufficiently large) such that, for every $R>R_{1}$, the following estimates hold:
		\begin{enumerate}[(1)]
			
			\item If $\Im(\lambda)=0$, then
			\begin{equation}
				\label{rel-l2-est-rss}
				\int_{R<d(o,x)<2R} |f(x)|^{2}\,dx \;\geq\; C\,R.
			\end{equation}
			
			\item If $1\le p<2$ and $\Im(\lambda)\neq 0$, then
			\begin{equation}
				\label{rel-lp-est-rss}
				\int_{R<d(o,x)<2R} |f(x)|^{p}\,dx 
				\;\geq\;
				C\begin{cases}
					e^{\,p(\gamma_{p}\rho - |\Im(\lambda)|)\,R}, 
					& \text{if } |\Im(\lambda)|\neq \gamma_{p}\rho, \\[6pt]
					R, & \text{if } |\Im(\lambda)|=\gamma_{p}\rho,
				\end{cases}
			\end{equation}
			where $\displaystyle \gamma_{p}=\frac{2}{p}-1$.
			
		\end{enumerate}
		
		\medskip
		\noindent
		Here $C>0$ denotes a constant depending on $X$, $\lambda$, $p$ (when relevant), and $f$, but independent of $R$.
	\end{thm}

	\begin{rem}
		\label{remark-rellich-asymptotic}
		The above theorem, in the spirit of Rellich’s result (Theorem~\ref{rel-org}), 
		clarifies the quantitative behavior of eigenfunctions across different spectral regimes. 
		In contrast to the Euclidean setting, the present framework naturally includes the case of non-real spectral parameters, 
		which arise naturally in the analysis on symmetric spaces of  noncompact type. 
		The relevance of such spectral regions becomes apparent when one examines the structure of the $L^p(X)$–spectrum of $\Delta_X$ (viz. \eqref{lp-spec}).
		
		When $\Im(\lambda) = 0$, the theorem shows that the $L^2$–mass of eigenfunctions in geodesic annuli grows at least {linearly} in $R$, 
		in agreement with the classical Euclidean behavior. 
		For $\Im(\lambda) \neq 0$ and $|\Im(\lambda)| < \gamma_p \rho$, the $L^p$–mass exhibits at least 
		{exponential growth} in $R$, reflecting the dominant influence of the exponential volume growth of $X$ in this spectral region. 
		At the critical spectral value $|\Im(\lambda)| = \gamma_p \rho$, this growth transitions from exponential to {linear}. 
		In the supercritical regime $|\Im(\lambda)| > \gamma_p \rho$, the estimate reveals an {exponential decay} 
		of the $L^p$–mass with respect to $R$, indicating that the decay rate of the eigenfunction outweighs 
		the exponential volume growth of the underlying space. 
		
		This trichotomy—exponential growth, linear growth, and exponential decay—is intrinsic to spaces of exponential volume growth 
		and has no direct analogue in the Euclidean setting.
	\end{rem}
	
	Furthermore, when $|\Im(\lambda)| \leq \gamma_p \rho$, 
	the growing lower bounds for the integrals in the preceding theorem demonstrate the \emph{non-integrability} 
	of the corresponding eigenfunctions at infinity. 
	As a consequence, we obtain the following result, which constitutes our second main theorem.
	
	\begin{thm}
		\label{rel-main-thm}
		Let $1 \leq p \leq 2$, and let $\Omega$ be as in Theorem~\ref{rellich-asymptotic-rss}. 
		Let $\lambda \in \mathbb{C}\setminus i\mathbb{Z}$ satisfy $|\Im(\lambda)| \leq \gamma_p \rho$. 
		Assume that $f$ is a solution of \eqref{helmrss-intro} in $\Omega$ such that
		\begin{equation}
			\label{rel-X-Lp-est}
			\int_{\Omega} |f(x)|^p \, dx < \infty.
		\end{equation}
		Then $f = 0$ in $\Omega$. 
	\end{thm}
This \textit{Rellich-type uniqueness} result may be viewed as a Liouville-type theorem for eigenfunctions in domains with a cavity. In light of the structure of the \(L^{p}\)-spectrum~\eqref{lp-spec}, it implies, in particular, the absence of \(L^{p}\)-point spectrum for the Laplace--Beltrami operator on such exterior domains~\(\Omega\).

\medskip
\noindent
\begin{rem}
	Several remarks are in order.
	\begin{enumerate}
		
		\item \textbf{Sharpness.}
		The theorem is sharp. For \(p > 2\) and any \(\lambda \in \mathbb{C} \setminus i\mathbb{Z}\), there exist eigenfunctions in \(\Omega\) with eigenvalue \(-(\lambda^{2}+\rho^{2})\) satisfying~\eqref{rel-X-Lp-est}. Likewise, for \(1 \le p \le 2\) and \(|\Im(\lambda)| > \gamma_p \rho\), one can construct eigenfunctions in \(\Omega\) satisfying~\eqref{rel-X-Lp-est}. This confirms the sharpness of the spectral region in the theorem (see Remark~\ref{rem:sharpness-lp-est}).
		
		In establishing this, we employ the elementary spherical function~$\varphi_{\lambda}$ 
		and the function~$\Phi_{\lambda}$ appearing in the Harish–Chandra expansion of~$\varphi_{\lambda}$ 
		(see~\eqref{HC-expan}).  
		We recall that $\Phi_{\lambda}$ is an eigenfunction of~$\Delta_X$, defined away from the identity, 
		with the same eigenvalue as~$\varphi_{\lambda}$ (See Section~\ref{sec:eigenfunction} for details).  
		In fact, in the analysis of eigenfunctions on~$\Omega$, $\Phi_{\lambda}$ serves as the principal model, 
		while in the full space~$X$, the role is played by~$\varphi_{\lambda}$.
		\item \textbf{On the restriction \(\lambda \notin i\mathbb{Z}\).}
		Analytically speaking, our approach in proving Theorem~\ref{rellich-asymptotic-rss} 
		is based on reducing the Helmholtz equation~\eqref{helmrss-intro} to a hypergeometric 
		differential equation with a regular singularity, via spherical harmonic expansion. 
		The exclusion in the spectral parameter $\lambda$ is then necessary to ensure the 
		existence and linear independence of the two fundamental solutions near the regular 
		singular point at $1$; see Remark~\ref{rem:exclusion-dis} for further details.

		From a harmonic-analytic perspective, this reflects the structure of the Harish--Chandra expansion~\eqref{HC-expan}: the \(c\)-function $c(\pm \lambda)$ is meromorphic with poles along \(i\mathbb{Z}\), and the independence of the components \(\Phi_{\pm\lambda}\) fails at these singular parameters. Moreover, the Poisson transform is not injective for certain \(\lambda \in i\mathbb{Z}\). Thus, the condition \(\lambda \notin i\mathbb{Z}\) is intrinsic to the spectral theory and ensures both analytic regularity and nondegeneracy (see Section~\ref{sec:eigenfunction}).
		
		\item \textbf{Critical exponent.}
		In rank-one symmetric spaces of noncompact type, the critical exponent is \(p=2\), independent of the dimension of~\(X\). This contrasts with the Euclidean case \(\mathbb{R}^{n}\), where the critical exponent is \(p=\tfrac{2n}{n-1}>2\) (see Theorem~\ref{BG1-thm}). This difference reflects the role of exponential volume growth and the asymptotics of spherical functions.
		
	\end{enumerate}
\end{rem}
	We next compare our result with two related works on Rellich-type theorems in geometric settings.
	
	Recently, Ballmann--Mukherjee--Polymerakis~\cite{BMP}, among other things, established an \(L^2\)-Rellich type uniqueness theorem for complements of horoballs in negatively curved asymptotically harmonic Hadamard manifolds, using (Donnelly–Garofalo) Rellich-type identities, integration by parts, and positivity properties of Busemann functions together with self-adjoint elliptic boundary conditions; see \cite[Theorem 3.2]{BMP}. Their result is formulated in terms of the absence of \(L^2\)-point spectrum for the Laplacian on horoball complements with self-adjoint elliptic boundary condition. The geometric setting considered in \cite{BMP} includes rank one Riemannian symmetric spaces of noncompact type. However, horoballs, defined as sublevel sets of Busemann functions, are typically noncompact; in rank-one symmetric spaces \(G/K\), they are bounded by horospheres arising as orbits of conjugates of \(N\) in the Iwasawa decomposition \(G=KAN\) (see, e.g., \cite[Ch.~II]{H2}). In contrast, we consider complements of geodesic balls centered at the identity coset \(o=eK\). Although an eigenfunction on the complement of a geodesic ball restricts to an eigenfunction on the complement of a horoball, the restricted function need not satisfy the boundary conditions required in the argument of \cite{BMP}; consequently, our \(L^2\)-uniqueness theorem does not follow directly from their work. Moreover, the analytical approaches are fundamentally different. While the method in \cite{BMP} is geometric and based on Rellich-type identity (See for instance \cite[Eq. (3.1)]{BMP}), our analysis relies on spherical harmonic expansions, reduction to hypergeometric differential equations, and asymptotic analysis of radial solutions. This yields the principal novelty of the present work: sharp quantitative lower bounds for the annular \(L^p\)-mass of eigenfunctions, leading to Rellich-type uniqueness results on the exterior of geodesic balls without imposing any boundary conditions.

	Another related result is due to Chen and Liu~\cite{CL}, who developed a scattering theoretic framework on real hyperbolic spaces \(\mathbb H^{n}\) and established, among other things, an \(L^{2}\)-Rellich type uniqueness theorem for solutions of the Helmholtz equation in the exterior of geodesic balls. The result of \cite{CL} is confined to the \(L^{2}\)-setting on real hyperbolic spaces and requires a stronger decay assumption at infinity. In contrast, our results apply to all rank one Riemannian symmetric spaces of noncompact type and yield sharp quantitative lower bounds for the annular \(L^{p}\)-mass of eigenfunctions, leading to Rellich-type uniqueness theorems for a significantly broader range of \(L^{p}\) spaces.
	
	\medskip
	The above discussion does not exhaust the Rellich-type phenomenon in the context of symmetric spaces. We next show that the Rellich-type uniqueness result admits a significant strengthening for non-real spectral parameters. More precisely, our third main result is as follows.
	
	\begin{thm}
		\label{rel-rss-weak-lp}
		Fix \(1 \leq p < 2\), and let \(\Omega\) be as in Theorem~\ref{rellich-asymptotic-rss}. 
		Let \(f\) be a solution of~\eqref{helmrss-intro} in \(\Omega\) with \(\lambda \in \mathbb{C} \setminus i\mathbb{Z}\). Then:
		
		\begin{enumerate}
			\item If \(\lambda = \alpha \pm i \gamma_p \rho\) with \(\alpha \in \mathbb{R}\), then \(f=0\) whenever \(f \in L^q(\Omega)\) for some \(1 \leq q \leq p\).
			
			\item If \(\lambda = \alpha \pm i \gamma_p \rho\) with \(\alpha \in \mathbb{R}\), then \(f=0\) whenever \(f \in L^{q,\infty}(\Omega)\) for some \(1 < q < p\).
			
			\item If \(\lambda = \alpha \pm i \gamma_q \rho\) with \(\alpha \in \mathbb{R}\) and \(p < q <2\), then \(f=0\) whenever \(f \in L^{p,\infty}(\Omega)\).
		\end{enumerate}
		In particular, if $1<p<2$, \(0\neq |\Im(\lambda)|<\gamma_p\rho\), and $f\in L^{p,\infty}(\Omega)$, then $f=0$ in $\Omega.$
	\end{thm}
	
	\noindent
	In particular, uniqueness persists up to the boundary of the \(L^{p}\)-spectrum: even at critical spectral parameters, weak integrability already forces triviality.
	
	\begin{rem}
		We record a few remarks.
		
		\begin{enumerate}
			\item The result is sharp; see Remark~\ref{weak-lp-sharpness} for details.
			
			\item Related weak \(L^{p}\)-estimates for eigenfunctions on the full space were obtained in Kumar--Ray--Sarkar~\cite[Proposition~3.1.1]{KRS}. Their approach relies on mean value type arguments (see \cite{H2}), which do not appear to extend to the exterior domain \(\Omega\), making the present setting more delicate.
			
			\item For comparison with the \(L^{p}\) Rellich-type uniqueness result (Theorem \ref{rel-main-thm}), observe that when \(\lambda\) lies in the interior of the \(p\)-strip (i.e., the interior of the \(L^{p}\)-spectrum), there are no nontrivial eigenfunctions even in weak \(L^{p}\).
		\end{enumerate}
	\end{rem}
	
	\medskip 
	Next, we turn to the latter part of this article, where we record some interesting phenomena exhibited by eigenfunctions in domains containing a cavity. 
	Setting aside the feature of the $L^p$ spectrum discussed above, Theorems~\ref{rel-main-thm}, and \ref{rel-rss-weak-lp} may be viewed as  standalone results describing the admissible size estimates of eigenfunctions in~$\Omega$ in terms of their $L^p$ and weak $L^p$ norms. 
	There are, nevertheless, several well-known results in the literature that characterize eigenfunctions of~$\Delta_X$ in $X$ satisfying various size estimates as Poisson transforms of suitable boundary data---for instance, functions belonging to appropriate Lebesgue classes or finite measures---on the boundary $K/M$ of~$X$. 
	In what follows, we focus on two such results involving Hardy-type norm estimates, which we now state. 
	\begin{thm}[\cite{BOS,Bou-Sami,Ion-Pois-1,KRS}]
		\label{poisson-charac}
		Let $\lambda\in \mathbb{C}$, and $u\in C^\infty(X)$ be such that $\Delta_Xu+(\lambda^2+\rho^2)u=0$ in $X.$ 
		\begin{enumerate}
			\item  Let $\Im(\lambda)<0$ or $\lambda=0$, and $1<p\leq \infty.$ Then there exists $F\in L^p(K/M)$ such that $u=\mathcal{P}_\la F$ if and only if 
			\begin{equation}
				\sup_{t>0}~\frac{1}{\varphi_{i\Im(\lambda)}(a_t)}\left(\int_K |u(ka_t)|^p~dk\right)^{\frac1p} <\infty. 
			\end{equation}
			If $p=1$, then $u=\mathcal{P}_{\la}\mu$ for some signed measure $\mu$ on $K/M.$
			\item Let $\Im(\lambda)=0$ with $\lambda\neq 0$, and $1\leq p\leq \infty.$ Then there exists $F\in L^p(K/M)$ such that $u=\mathcal{P}_{\la}F$ if and only if 
			$$\sup_{t>0}~e^{\rho t}\left(\int_K |u(ka_t)|^p~dk\right)^{\frac1p} <\infty.$$
		\end{enumerate}
		Here $\mathcal{P}_\la F$ stands for the Poisson transform on $F$, see \eqref{poisson-def}.
	\end{thm}
	The part (1) of the above theorem was established by Ben Said et al. \cite{BOS}, while Bousseria-Sami \cite{Bou-Sami} and Ionescu \cite{Ion-Pois-1} proved the part (2) for the range $p\geq 2$. It was further conjectured in \cite{Bou-Sami} that the result should also hold for $1<p<2$, which was subsequently confirmed affirmatively by Kumar-Ray-Sarkar\cite{KRS}.
	
	Motivated by the role of such Hardy-type norms in the characterization of eigenfunctions, we establish—as a byproduct of our study of Rellich-type theorems in the present setting—the following result concerning admissible size estimates for eigenfunctions in the complement of a cavity.
	\begin{thm}
		\label{rel-hardy-norm}
		Let $\lambda \in \mathbb{C}\setminus i\mathbb{Z}$, $\epsilon > 0$, $1\leq p\leq \infty$ and let $\Omega$ be as in Theorem~\ref{rellich-asymptotic-rss}. 
		Assume that $f \in C^2(\Omega)$ is a solution of \eqref{helmrss-intro} in $\Omega$. Then  there exists an explicit positive function $\psi_{\lambda}$ (See \eqref{def-psi}) such that if 
		\begin{equation}
			\label{rel-hardy-norm-est}
			\sup_{t > R_0} 
			\frac{t^{\epsilon}}{\psi_{\lambda}(a_t)} 
			\left( \int_K |f(k a_t)|^p \, dk \right)^{\frac{1}{p}} < \infty,
		\end{equation}
		(with the usual modification for $p=\infty$) then $f = 0$ in $\Omega$.
		Moreover, the condition $\epsilon > 0$ is \emph{sharp} in the sense that if $\epsilon = 0$, there exists a nontrivial eigenfunction $f$ in $\Omega$ satisfying the  estimate \eqref{rel-hardy-norm-est}
	\end{thm}
	Note that while Theorem~\ref{poisson-charac} provides an $L^p$–characterization of eigenfunctions on the entire symmetric space $X$ via the Poisson transform—where the growth of $u$ along geodesic spheres is governed by Hardy-type norms—our result, Theorem~\ref{rel-hardy-norm}, concerns eigenfunctions on the \emph{complement of a compact cavity} in $X$ and establishes a \emph{Rellich-type uniqueness} theorem. We also highlight an additional structural feature, stated in the following remark.
	\begin{rem}
		The function $\psi_{\lambda}$ in our result satisfies 
		$\psi_{\lambda}(a_t) \asymp |\Phi_{-i\,\Im(\lambda)}(a_t)|$ for large $t$ whenever $\Im(\lambda)<0$ 
		(see~\eqref{psi-asymp}), 
		where $\Phi_{-i\,\Im(\lambda)}$ denotes the term appearing in the Harish–Chandra expansion 
		of $\varphi_{\,i\,\Im(\lambda)}$ (see~\eqref{HC-expan}) that was used in Theorem~\ref{poisson-charac} 
		to characterize eigenfunctions. 
		This function is itself an eigenfunction in $\Omega$ (though not necessarily on $X$) 
		corresponding to the same eigenvalue as $\varphi_{\,i\,\Im(\lambda)}$; 
		see Section~\ref{sec:eigenfunction} for details. 
		Moreover, when $\Im(\lambda)=0$, one has $\psi_{\lambda}(a_t)=e^{-\rho t}$ 
		(see~\eqref{def-psi}). 
		These observations highlight a natural correspondence between our result and the 
		Poisson-type characterization discussed earlier (Theorem~\ref{poisson-charac}), albeit arising in a distinct geometric setting.
	\end{rem}
	We conclude the introduction with a brief outline of the paper. 
	Section~2 collects the necessary preliminaries on rank-one Riemannian symmetric spaces of noncompact type. 
	In Section~3, we present the proofs of the main results. 
	The paper ends with some concluding remarks, highlighting related questions and possible directions for future research in Section~4.
	
	\subsection*{Notations} We employ the notation $C,c$, $C_1,c_1$, $C_2, c_2$, {etc.}, for  constants whose values may differ at each occurrence. For two nonnegative functions $f$ and $g$, we write 
	$f \lesssim g$ (respectively, $f \gtrsim g$) to mean that there exists a constant 
	$C>0$ such that $f \le Cg$ (respectively, $f \ge Cg$).  
	The notation $f \asymp g$ signifies that both $f \lesssim g$ and $f \gtrsim g$ hold.
	We also use standard asymptotic notations:  
	$f \sim g$ means that 
	$
	\lim_{t\to\infty} \frac{f(t)}{g(t)} = 1,
	$
	while $f = o(g)$ means that 
	$
	\lim_{t\to\infty} \frac{f(t)}{g(t)} = 0.
	$
	
	\section{Riemannian symmetric spaces of  noncompact type}\label{sec:rss}
	In this section, we introduce the notations and relevant concepts related to the theory of Riemannian symmetric spaces of  noncompact type. Most of these are standard and widely available in sources such as \cite{H1,H2}.
	
	Let $X$ be a Riemannian symmetric space of  noncompact type. Then $X$ can be realized as the quotient space $X=G/K$ where $G$ is a connected,  noncompact real semi-simple Lie group with finite center and $K$ a maximal compact subgroup of $G.$  The group $G$  naturally acts on $X$ by left translations $l_g: xK\mapsto gxK,~ g\in G.$ We denote the base point $eK$ simply by $o.$ 
	
	\subsection{Basic structure theory and Riemannian metric} We denote the Lie algebras of $G$ and $K$ by the corresponding Gothic fraktur notations $\gf$, and $\kf$, respectively. The corresponding Cartan decomposition reads as $\gf=\kf\bigoplus\pf$ where $\pf$ is the orthogonal complement of $\kf$ with respect to the Cartan-Killing form $B$ of $\gf.$  Let $\af$ be a maximal abelian subspace of $\pf.$ The dimension of $\af$ which is independent of its choice is called the rank of the symmetric space $X.$ In what follows, we assume that $\dim \af =1.$ Now denoting the real dual of $\af$ by $\af^*$, we let $\Sigma\subset \af^*$ be the restricted system of roots of the pair $(\gf, \af).$ For $\alpha\in \Sigma$, let $\gf_\alpha$ denote the associated root space and we set $\dim \gf_\alpha=m_{\alpha}.$ Fix a positive Weyl chamber $\af_+\subset \af$ which leads to the corresponding set of positive roots  $\Sigma^+$. Let $\mathfrak n=\sum_{\alpha\in \Sigma^+} \gf_\alpha$. Then  $N=\exp \mathfrak n$ is a simply connected nilpotent Lie group and  $A=\exp \af$ is an abelian group. Let $M$ be the centralizer of $A$ in $K$. Now the usual product map $(k, a, n)\mapsto kan$ from $K\times A\times N$ to $G$ is an analytic diffeomorphism which is not a group homomorphism.  This leads to the  Iwasawa decomposition: $G=KAN$. In view of  this decomposition an element $g\in G$ can be uniquely written as $g=k \exp(H(g)) n$ where $k\in K, n\in N$ and $H(g)\in \af$.
	The group $G$ also admits a polar decomposition $G=K \overline{A^+} K$ a.k.a. the Cartan decomposition where $A^+:=\exp \af_+.$   Let $M$ be the centralizer of $A$
	in $K$. Then we note that the groups $M$ and $A$ normalize $N$ and $ K/M$ is the Furstenberg boundary of $X.$
	
	Now since in our case $\dim \af=1$, it is well-known that $\Sigma$ is either $\{\pm \gamma\}$ or $\{\pm \gamma, \pm 2\gamma\}$ where $\gamma$ is a positive root, and the associated Weyl group $W=\{\pm I\}.$ It is customary to denote the half sum of the positive roots (counted with multiplicities) by $\rho$, i.e., $\rho=\frac12 (m_{\gamma}+2m_{2\gamma})\gamma.$ Let $H_0\in \af$ denote the unique element so that $\gamma (H_0)=1$ leading to the identification of $\af$ with $\mathbb{R}$ via the map $t\mapsto t H_0~(t\in\mathbb{R}).$ Consequently, the maps $t\mapsto t\gamma~(t\in\mathbb{R})$, and $z\mapsto z\gamma~(z\in\mathbb{C})$ identifies $\af^*$ with $\mathbb{R}$ and its complexification $\af^*_{\mathbb{C}}$ with $\mathbb{C}$, respectively. 
	Under this identification, we follow the customary abuse of notation and denote
	\begin{equation}\label{rho-def}
		\rho = \rho(H_0) = \tfrac{1}{2}(m_{\gamma}+2m_{2\gamma}).
	\end{equation}
	Moreover, the group $A$ can be realized as $A=\{a_t:=\exp (t H_0): t\in\mathbb{R}\}$ and subsequently $A^+=\{a_t: t>0\}.$ The Lebesgue measure on $\mathbb{R}$ induces a Haar measure on $A$ which we simply denote by $da_t=dt.$
	
	In view of the Iwasawa decomposition, the map 
	\[
	A \times N \longrightarrow X, \qquad (a_t,n) \mapsto a_tn \cdot o,
	\]
	is a surjective diffeomorphism. Consequently, the dimension of $X$ is given by 
	$
	n = \dim(X) = m_{\gamma} + m_{2\gamma} + 1.
	$
	Moreover, as is well known from the Cartan decomposition, the double cosets $Ka_tK$ and $Ka_sK$ coincide if and only if $t = \pm s$. Hence, every point $x \in X$ can be written in the form 
	$
	x = k a_t \cdot o,
	$
	for some $k \in K$ and a unique $t \ge 0$.
	
	The Killing form $B$ on $\gf$ induces a $K$-invariant inner product on $\pf$, and thus defines a $G$-invariant Riemannian metric on $X$ together with an associated Riemannian measure, denoted by $dx$. Viewing a function $f$ on $X$ as a right $K$-invariant function on $G$, the Haar measure $dg$ on $G$ can be normalized so that
	\[
	\int_X f(x)\,dx = \int_G f(g)\,dg.
	\]
	For any integrable function $f$ on $X$, the integration formula corresponding to the polar decomposition takes the form
	\begin{equation}
		\label{polar-int}
		\int_X f(x)\,dx = \int_K \int_0^{\infty} f(k a_t)\, J(t)\,dt\,dk,
	\end{equation}
	where $dk$ denotes the normalized Haar measure on $K$, and the Jacobian of the polar decomposition is given by
	\begin{equation}
		\label{def-J}
		J(t) = c_n (2\sinh t)^{m_{\gamma}+m_{2\gamma}} (\cosh t)^{m_{2\gamma}}.
	\end{equation}
	It then follows that for any $c>0$,
	\begin{equation}
		\label{est:J}
		J(t)\asymp e^{2\rho t} \quad\quad (t>c).
	\end{equation}
	The metric induced by the Killing form endows $X$ with the structure of a Riemannian manifold. The corresponding Riemannian distance is denoted by $d(\cdot,\cdot)$, and $\Delta_X$ denotes the associated Laplace–Beltrami operator. We now describe the explicit form of $\Delta_X$ in polar co-ordinates, which will play a central role in the subsequent analysis. For this, we closely follow Helgason\cite[p.309, Ch.II]{H1}.
	
	\subsection{Polar co-ordinates and Laplace-Beltrami operator}  
	For $\alpha\in \Sigma^+$, let 
	$$\kf_{\alpha}:=\{T\in \kf: ~(\ad(H))^2T=\alpha(H)^2T~~\text{for}~~H\in\af\},$$ and $\{T^{\alpha}_1,\hdots, T^{\alpha}_{m_{\alpha}}\}$ a basis of $\kf_{\alpha}$, orthonormal with respect to the negative of the Killing of $\gf.$ We set 
	\begin{equation}
		\label{defomega}
		\omega_{\alpha}:=\sum_{i=1}^{m_{\alpha}} T_i^{\alpha}\cdot T_i^{\alpha}
	\end{equation}
	where  the $T_i^{\alpha}$ are viewed as left-invariant differential operators on $G.$ Then it is easy to see that $\omega_{\alpha}$ is bi-invariant under $M.$
	As proved in \cite[Theorem 5.24]{H1}, the Laplace-Beltrami operator admits the  following representation in polar co-ordinates: 
	\begin{align}
		\label{LBexp}
		\Delta_{X}f(ka\cdot o)= &\left(\Delta_A+\sum_{\alpha\in \Sigma^+}m_{\alpha}(\coth \alpha) A_{\alpha} \right)_af(ka\cdot o) \nonumber \\
		& +\sum_{\alpha\in \Sigma^+} \sinh^{-2}(\alpha(\log a)) ((\Ad(a^{-1})\omega_{\alpha})f)(ka)
	\end{align}
	where $\Delta_A$ denotes the Laplacian on $A\cdot o,$ and $A_\alpha$ is given by $$\langle A_\alpha, H\rangle=\alpha(H)\quad\quad (H\in\mathfrak{a},~\alpha\in \Sigma^+).$$
	Also, $\Ad$ stand for the adjoint representation of $G$ on its Lie algebra $\mathfrak{g}$, and  for any $X\in \mathfrak{g}$, we  have
	$$\exp(\Ad(g))X=g\exp(X) g^{-1},\ ~~~\Ad(\exp(X))=e^{\ad X}\quad\quad (g\in G),$$
	where $\ad$ is the adjoint representation of $\mathfrak{g}$ defined by $\ad X(Y)=[X,Y]$ for $X, Y\in \mathfrak{g}.$

	Next, we describe a useful closed form expression for $\Delta_{X}$ in the co-ordinate system $\exp(tH_0)\rightarrow t.$ First note that, since for $f\in C^\infty(A),$ $$H_0f(\exp(tH_0)=\frac{d}{dt}f(\exp(tH_0),$$ $H_0$ is identified with $\frac{d}{dt}$ in this co-ordinate system. Now observe that, as $\mathfrak{a}$ is one-dimensional, there exists $C_\gamma$ such that $A_\gamma=C_\gamma H_0.$ Now observe that $$B(H_0, H_0)=2(m_\gamma \gamma(H_0)^2+m_{2\gamma}(2\gamma(H_0))^2)=2(m_\gamma+4m_{2\gamma})$$
	whence $A_\gamma=(2(m_\gamma+4m_{2\gamma}))^{-1}~H_0=(2(m_\gamma+4m_{2\gamma}))^{-1}~\frac{d}{dt}.$ Similarly, $A_{2\gamma}= (m_\gamma+4m_{2\gamma})^{-1}~\frac{d}{dt}.$ Subsequently, the Laplacian $\Delta_A$ on $A$ takes the form 
	$(2(m_\gamma+4m_{2\gamma}))^{-1}~\frac{d^2}{dt^2}.$ As a result, the Laplace-Beltrami \eqref{LBexp} takes the following form
	\begin{align}
		\label{explicitLB-rk1}
		\Delta_X= \frac12(m_\gamma+4m_{2\gamma})^{-1}&\left( \frac{d^2}{d t^2} + (m_{\gamma} \coth t+ 2m_{\gamma} \coth (2t))\frac{d}{d t}\right)\nonumber \\ &+\sinh^{-2}(t)(\Ad(a_t^{-1})\omega_{\gamma})+\sinh^{-2}(2t)(\Ad(a_t^{-1})\omega_{2\gamma}).
	\end{align}
	\subsection{Spherical harmonics}
	We now describe the system of spherical harmonics adapted to the present setting. 
	To this end, we make use of certain irreducible representations of $K$ admitting $M$–fixed vectors. 
	Let $\widehat{K_0}$ denote the set of all irreducible unitary representations of $K$ having nontrivial $M$–fixed vectors. 
	For $\delta \in \widehat{K_0}$, let $V_{\delta}$ be the finite-dimensional Hilbert space on which $\delta$ is realized. 
	As shown in~\cite[Theorem~6]{Ks}, the subspace of $V_{\delta}$ consisting of $M$–fixed vectors is one-dimensional. 
	Hence, $V_{\delta}$ contains a unique normalized $M$–fixed vector, denoted by~$v_1$. 
	Let $\{v_1,v_2,\dots,v_{d_{\delta}}\}$ be an orthonormal basis of $V_{\delta}$. 
	For each $\delta \in \widehat{K_0}$ and $1 \le j \le d_{\delta}$, define
	\begin{equation}
		\label{Sph-Harm-def}
		Y_{\delta,j}(kM) = (v_j,\, \delta(k)v_1)_{L^2(K/M)}, \qquad kM \in K/M.
	\end{equation}
	It follows immediately that $Y_{\delta,1}(eK)=1$ and that $Y_{\delta,1}$ is $M$–invariant. 
	
	\begin{prop}
		\label{H2}
		The collection $\{Y_{\delta,j} : 1 \le j \le d_{\delta},\, \delta \in \widehat{K_0}\}$ forms an orthonormal basis of $L^2(K/M)$.
	\end{prop}
	
	\noindent
	For a proof, we refer the reader to~\cite[Theorem~3.5, Section~3, Chapter~5]{H2}; see also the discussion at the beginning of~\cite[Section~2, Chapter~3]{H1}.
	
	\smallskip
	An explicit realization of $\widehat{K_0}$ can be obtained by identifying $K/M$ with the unit sphere in~$\mathfrak{p}$. 
	Let $\mathcal{H}^m$ denote the space of homogeneous harmonic polynomials of degree~$m$ restricted to the unit sphere. 
	Then the classical spherical harmonic decomposition takes the form
	\[
	L^2(K/M) = \bigoplus_{m=0}^{\infty} \mathcal{H}^m.
	\]
	Each $V_{\delta}$ is contained in some $\mathcal{H}^m$, and the functions $Y_{\delta,j}$ may thus be identified with spherical harmonics. 
	
	\begin{thm}[{\cite[Theorem~11.2]{H2}}]
		\label{omega-eigenvalue}
		For each $\delta \in \widehat{K_0}$, there exist unique integers $q \ge p \ge 0$ such that
		\begin{align*}
			\delta(\omega_{\gamma})\big|_{V^M_{\delta}}
			&= \frac{p(p+m_{2\gamma}-1) - q(q+m_{\gamma}+m_{2\gamma}-1)}{2(m_{\gamma}+4m_{2\gamma})}, \\[4pt]
			\delta(\omega_{2\gamma})\big|_{V^M_{\delta}}
			&= -\,\frac{2p(p+m_{2\gamma}-1)}{m_{\gamma}+4m_{2\gamma}},
		\end{align*}
		with the convention that $p=0$ whenever $m_{2\gamma}=0$.
	\end{thm}
	
	\subsection{Eigenfunctions of the Laplace-Beltrami operator $\Delta_X$}\label{sec:eigenfunction} In this subsection, we briefly recall certain eigenfunctions whose estimates will play a central role in our analysis. We begin with the \emph{elementary spherical functions}. For $\lambda \in \mathbb{C}$, these are defined by  
	\[
	\varphi_{\lambda}(x) = \int_{K} e^{-(i\lambda + \rho) H(x^{-1}k)}\, dk, \qquad x \in X.
	\]
	It is well known that $\varphi_{\lambda}$ is a $K$-biinvariant function on $G$ satisfying  
	\[
	\varphi_{\lambda}(o) = 1,~ \varphi_{\lambda} = \varphi_{-\lambda} \quad \text{for all } \lambda \in \mathbb{C}, 
	\qquad \text{and} \qquad 
	\varphi_{\lambda}(x) = \varphi_{\lambda}(x^{-1}) \quad \text{for all } x \in X.
	\]
	Moreover, $\varphi_{\lambda}$ is an eigenfunction of the Laplace--Beltrami operator $\Delta_X$ with eigenvalue $-(\lambda^2 + \rho^2)$, that is,
	$
	\Delta_X \varphi_{\lambda} + (\lambda^2 + \rho^2)\varphi_{\lambda} = 0.
	$
	It follows immediately from the integral representation that the pointwise estimate $
	|\varphi_{\lambda}| \leq \varphi_{i\,\Im(\lambda)} $
	holds for all $\lambda \in \mathbb{C}$.
	
	To obtain more precise information regarding the behavior of $\varphi_{\lambda}$ away from the identity, we recall the classical \emph{Harish-Chandra expansion} (see Stanton--Tomas~\cite[Theorem~3.1]{ST}; see also \cite[p. 7--8]{K2}):
	\begin{equation}\label{HC-expan}
		\varphi_{\lambda}(a_t) = c(\lambda)\Phi_{\lambda}(a_t) + c(-\lambda)\Phi_{-\lambda}(a_t), \quad\lambda\in \mathbb{C}\setminus i\mathbb{Z},
	\end{equation}
	where
	\begin{equation}\label{Phi-la}
		\Phi_{\lambda}(a_t) = e^{(i\lambda - \rho)t} \sum_{k=0}^{\infty} \Gamma_k(\lambda) e^{-2kt}.
	\end{equation}
	The coefficients $\Gamma_k(\lambda)$ in the above series are defined recursively by $\Gamma_0(\lambda) \equiv 1$ , and for $k\geq 1$
	\[ 
	(k+1)(k+1 - i\lambda)\Gamma_{k+1} = 
	\sum_{j=0}^{k} \tfrac{m_{\gamma}}{2}(\rho + 2j - i\lambda)\Gamma_j 
	+ \sum_{\substack{j = k + 1 - 2l \\ j,l \ge 0}} m_{2\gamma}(\rho + 2j - i\lambda)\Gamma_j.
	\]
	Here $c(\lambda)$ denotes the \emph{Harish--Chandra $c$–function}, given by 
	\cite[(3.2)]{ST},
	\begin{equation}
		\label{HC-c-fns}
		c(\lambda)
		= \Gamma\!\left(\frac{m_{\gamma}}{2}\right)
		\Gamma\!\left(\frac{m_{2\gamma}}{2}\right)
		\frac{%
			\Gamma(i\lambda)\,
			\Gamma\!\left(\frac{m_{\gamma}+i\lambda}{2}\right)}%
		{%
			\Gamma\!\left(\frac{m_{\gamma}}{2}+i\lambda\right)\,
			\Gamma\!\left(\frac{\rho+i\lambda}{2}\right)}.
	\end{equation}
	
	We summarize below several key features of this expansion:
	
	\begin{enumerate}
		\item[\textbf{(1)}] \textbf{The restriction $\lambda\in \mathbb{C}\setminus i\mathbb{Z}$.}  
		The assumption $\lambda\notin i\mathbb{Z}$ ensures that all the coefficients 
		$\Gamma_k(\lambda)$ and $\Gamma_k(-\lambda)$ appearing in the definitions of 
		$\Phi_{\lambda}$ and $\Phi_{-\lambda}$ are well-defined and that the two 
		solutions remain linearly independent.
		
		From the expression~\eqref{HC-c-fns}, it is clear that 
		$c(\lambda)$ is meromorphic with possible poles arising from the 
		\emph{Gamma functions in the numerator} of \eqref{HC-c-fns}.
		Since $\Gamma(z)$ has poles precisely at $z = 0,-1,-2,\ldots$, 
		these occur exactly when $\lambda\in i\mathbb{Z}_{\geq 0}.$ Consequently, both $c(-\la)$, and $c(\la)$ have no poles in $\mathbb{C}\setminus i\mathbb{Z}.$
		
		\item[\textbf{(2)}] \textbf{The function $\Phi_{\lambda}$.}  
		The series~\eqref{Phi-la} defining $\Phi_{\lambda}$ converges for $t>c>0$, away from the identity 
		element $a_0 = e$.  
		Moreover, $\Phi_{\lambda}$ is an eigenfunction of $\Delta_X$ (away from the identity) 
		with eigenvalue $-(\lambda^2 + \rho^2)$.  
		For further details, see Helgason~\cite[Chapter~IV, §5]{H1}.
		
		\item[\textbf{(3)}] \textbf{Estimates on the coefficients.}  
		Gangolli~\cite{Gangolli} established that there exist positive constants $c$ and $d$ such that 
		\[
		|\Gamma_k(\lambda)| \leq c\, k^d \qquad \text{for all } k \geq 1,
		\]
		and furthermore, that these bounds are sharp.
	\end{enumerate}
	
	Observe that, from \eqref{Phi-la} and the previously mentioned estimates for the coefficients $\Gamma_k$, it follows immediately that  
	\begin{equation}
		\label{Phi-asymp}
		\Phi_{\lambda}(a_t) \sim e^{(i\lambda - \rho)t} \qquad (t \to \infty).
	\end{equation}
	Consequently, for sufficiently large $t$, one obtains the asymptotic estimate  
	\begin{equation}
		\label{Phi-la-est}
		|\Phi_{\lambda}(a_t)| \asymp e^{(-\Im(\lambda) - \rho)t}.
	\end{equation}
	
	It is well known that the Harish-Chandra $c$-function $c(\lambda)$ has neither zeros nor poles in the half-plane $\Im(\lambda) < 0$.  
	Within this region, by combining the Harish-Chandra expansion \eqref{HC-expan} with \eqref{Phi-asymp}, one deduces that  
	\[
	\lim_{t \to \infty} e^{-(i\lambda - \rho)t}\, \varphi_{\lambda}(a_t) = c(\lambda).
	\]
	(See \cite[Theorem~6.14]{H1} for further details.)  
	From this observation, it follows that for $\Im(\lambda) \neq 0$ there exists $t_{\lambda} > 0$ such that  
	\begin{equation}
		\label{phi-la-est}
		|\varphi_{\lambda}(a_t)| \asymp e^{(|\Im(\lambda)| - \rho)t} \qquad (t > t_{\lambda}).
	\end{equation}
	Also, for $\lambda = 0$, one has the classical estimates  
	\[
	e^{-\rho t} \lesssim \varphi_{0}(a_t) \lesssim (1 + t)e^{-\rho t}, \qquad t > 0.
	\]
	In view of estimate~\eqref{phi-la-est}, for $p>2$, we have
	\begin{equation}
		\label{phi-la-lp-est}
		\varphi_{\lambda}\in L^p(X)\qquad\text{whenever}\qquad
		|\Im(\lambda)| < \Bigl(1-\frac{2}{p}\Bigr)\rho .
	\end{equation}
	Consequently, for every $p>2$ each spectral parameter $\lambda$ satisfying the above inequality
	gives rise to a nontrivial $L^p$–eigenfunction of $-\Delta_X$. 
	Equivalently, every point of the interior of the $L^p(X)$–spectrum of $-\Delta_X$ (see \eqref{lp-spec}) 
	is contained in the point spectrum of $-\Delta_X$ on $L^p(X)$.
	
	\medskip
	
	We now turn to a broader class of eigenfunctions.  
	Recall that the complex powers of the Poisson kernel,
	\[
	x \longmapsto e^{-(i\lambda + \rho)H(x^{-1}k)}, \qquad \lambda \in \mathbb{C},\ k \in K,
	\]
	are eigenfunctions of the Laplace--Beltrami operator $\Delta_X$ corresponding to the eigenvalue $-(\lambda^2 + \rho^2)$.  
	For $\lambda \in \mathbb{C}$ and $F \in L^1(K/M)$, the \emph{Poisson transform} of $F$ is defined by
	\begin{equation}
		\label{poisson-def}
		\mathcal{P}_{\lambda}F(x) = \int_{K} F(k)\, e^{-(i\lambda + \rho)H(x^{-1}k)}\, dk, \qquad x \in X.
	\end{equation}
	It is immediate that $\mathcal{P}_{\lambda}F$ satisfies the eigenvalue equation
	\[
	\Delta_X (\mathcal{P}_{\lambda}F) + (\lambda^2 + \rho^2)\, \mathcal{P}_{\lambda}F = 0.
	\]
	Moreover, when $F \equiv 1$, the definition yields $\mathcal{P}_{\lambda}F = \varphi_{\lambda}$.  
	In fact, by the general theory of harmonic analysis on symmetric spaces, every eigenfunction of the Laplace--Beltrami operator can be realized as a Poisson transform of an appropriate function (or, more generally, of an analytic functional) on the boundary $K/M$.  
	In particular, for the Dirac measure $\delta_{k_0}$ at $k_0 \in K$, one has
	\[
	\mathcal{P}_{\lambda}\delta_{k_0}(x) = e^{-(i\lambda + \rho)H(x^{-1}k_0)}.
	\]
	
	We recall that a parameter $\lambda\in \mathbb{C}$ is called \emph{simple} if the map 
	$F \mapsto \mathcal{P}_{\lambda}F$ (the Poisson transform) is injective 
	(see~\cite[p.~248]{H2}).  
	It is well known (see~\cite[Theorem~6.5]{H2}) that $\lambda$ is simple if and only if 
	$e(\lambda)\neq 0$, where $1/e(\lambda)$ is the denominator of the 
	Harish--Chandra $c$–function~\eqref{HC-c-fns}. For nonsimple parameters, the Poisson transform acquires a nontrivial kernel and the Harish--Chandra expansion degenerates.

	Finally, we recall that eigenfunctions of $\Delta_X$ satisfying Hardy-type growth or integrability conditions admit a precise characterization in terms of Poisson transforms.
	Such representations, formulated via Hardy-type norm estimates, are of particular significance for our analysis and have already been stated in the Introduction (see Theorem~\ref{poisson-charac}).
	
	\section{Proofs of main results}
	
	In this section, we present the proofs of our main results, namely 
	Theorems~\ref{rellich-asymptotic-rss}, \ref{rel-main-thm}, and \ref{rel-hardy-norm}. 
	Before proceeding, we record a few notational remarks for the reader's convenience.
	
	Throughout this section, we denote by $\Omega$ the complement of a geodesic ball centered at $o$, defined by
	\begin{equation}
		\label{def-Omega}
		\Omega := \{ x \in X : d(o, x) > R_0 \}, \qquad \text{for some } R_0 > 0.
	\end{equation}
	Recall that functions on $X$ are regarded as right $K$–invariant functions on $G$. 
	Consequently, the expressions $f(ka_t \cdot o)$ and $f(ka_t)$ represent the same object; 
	hence, we shall use them interchangeably whenever no confusion arises. 
	Moreover, as is customary, we shall not distinguish between integration over $K$ and over $K/M$ for notational simplicity. Throughout the sequel, various constants will appear at different stages. 
	Unless explicitly stated otherwise, these constants may depend on structural parameters of the space $X$; 
	however, since such dependence plays no essential role in our analysis, 
	we shall omit it from the notation for brevity.
	
	We begin with the following proposition, which provides a fundamental estimate 
	for solutions of the Helmholtz-type equation in $\Omega$.
	
	\begin{prop}
		\label{main-prop}
		Let $\lambda \in \mathbb{C}\setminus i\mathbb{Z}$, and let $\Omega$ be as in~\eqref{def-Omega}. 
		Assume that $f \in C^2(\Omega)$ is a non-trivial solution of
		\begin{equation}
			\label{helmrss}
			\Delta_X f + (\lambda^2 + \rho^2) f = 0 \qquad \text{in } \Omega.
		\end{equation}
		Then, there exists $R_1>R_0$ such that, for every $1 \leq p \leq \infty$ (with the usual modification for $p = \infty$), the following estimates hold:
		\begin{enumerate}[(i)]
			\item \textbf{Case} $\boldsymbol{\Im(\lambda) \neq 0.}$
			There exists a constant $C > 0$, depending on $\lambda$ and $p$, such that
			\begin{equation}
				\label{est-la-non-0}
				\int_K |f(ka_t)|^p \, dk 
				\;\geq\;
				C \, e^{p(-| \Im(\lambda)| - \rho)t},
				\qquad (t > R_1).
			\end{equation}
			
			\item \textbf{Case} $\boldsymbol{\Im(\lambda) = 0}.$
			There exist constants $A, B, C > 0$, and $\theta\in \mathbb{R}$ depending on $\lambda$ and other constants related to the structure of $X$, such that
			\begin{equation}
				\label{est-la-0}
				\left( \int_K |f(ka_t)|^p \, dk \right)^{\tfrac{1}{p}}
				\geq
				C e^{-\rho t} 
				\big( A^2 + B^2 + 2AB \cos(2\lambda t + \theta) \big)^{1/2},\quad(t>R_0).
			\end{equation}
		\end{enumerate}
	\end{prop}
	
	\begin{proof}
		Let $f$ be as in the statement of the theorem. Using polar co-ordinates $X_{reg}=K A^+\cdot o$, there exists $t_0>R_0$ such that $f\neq 0$ on the sphere of radius $t_0$ viz. $S_{t_0}=\{k a_{t_0}\cdot o: k\in K\}.$ Therefore, in view of the Proposition \ref{H2}, it follows that there exists $\delta\in \widehat{K_0}$ and $1\leq j\leq d_{\delta}$ such that  the function 
		\begin{equation}
			\label{def-u}
			u(t):=\int_{K} f(ka_t\cdot o)Y_{\delta,j}(k) ~dk\quad\quad (t>R_0)
		\end{equation}
		is not identically zero as $u(t_0)\neq 0.$  Now in view of the hypothesis \eqref{helmrss}, we see that 
		\begin{align}
			\label{sph-coeff-Helm}
			0&=\int_{K} \bigg(\Delta_{X}f(ka_t\cdot o)+(\lambda^2+\rho^2)f(ka_t\cdot o)\bigg)Y_{\delta,j}(k) ~dk\nonumber\\
			&= \int_{K} \Delta_{X}f(ka_t\cdot o)Y_{\delta,j}(k) ~dk+(\lambda^2+\rho^2) \int_{K/M} f(ka_t\cdot o)Y_{\delta,j}(k) ~dk\nonumber\\
			&=  \int_{K} \Delta_{X}f(ka_t\cdot o)Y_{\delta,j}(k) ~dk+(\lambda^2+\rho^2)~u(t).
		\end{align}
		But the explicit form of the Laplace-Beltrami operator \eqref{explicitLB-rk1} then shows that 
		\begin{align}
			\label{sph-coeff-laplacef}
			&\nonumber\int_{K} \Delta_{X}f(ka_t\cdot o)Y_{\delta,j}(k) ~dk\\
			&= \frac12(m_\gamma+4m_{2\gamma})^{-1}\int_{K} \left( \frac{d^2}{d t^2} + (m_{\gamma} \coth t+ 2m_{\gamma} \coth (2t))\frac{d}{d t}\right)f(ka_t\cdot o) Y_{\delta,j}(k) dk\\
			&\nonumber\hspace*{1cm} + \int_{K}\left(\sinh^{-2}(t)((\Ad(a_t^{-1})\omega_{\gamma})f)(ka_t)+\sinh^{-2}(2t)((\Ad(a_t^{-1})\omega_{2\gamma})f)(ka_t)\right) Y_{\delta,j}(k)~dk
		\end{align}
		In order to simplify the last expression, for any $\alpha\in \{\gamma,2\gamma\}$, recalling the definition of $\omega_{\alpha}$ (See \eqref{defomega}), and using the linearity of $\Ad$, we notice that
		\begin{align*}
			&\int_K ((\Ad(a_t^{-1})\omega_{\alpha})f)(ka_t) Y_{\delta,j}(k)~dk\\
			&=\sum_{i=1}^{m_{\alpha}} \int_K ((\Ad(a_t^{-1})(T_i^\alpha)^2)f)(ka_t) Y_{\delta,j}(k)~dk\\
			&=\sum_{i=1}^{m_{\alpha}} \int_K \left( \frac{d^2}{ds^2}\Bigg|_{s=0}f(ka_t\cdot \exp(s\Ad(a_t^{-1})T^\alpha_i)\cdot o)\right) Y_{\delta, j}(k)~dk 
		\end{align*}
		which, upon using $$\exp(s\Ad(a_t^{-1})T^\alpha_i)=a_t^{-1} \exp(sT_i^\alpha)a_t,$$
		transforms to
		\begin{align*}
			&\sum_{i=1}^{m_{\alpha}} \int_K \left( \frac{d^2}{ds^2}\Bigg|_{s=0}f(k\cdot \exp(sT_i^\alpha)a_t\cdot o)\right) Y_{\delta, j}(k)~dk\\
			&= \sum_{i=1}^{m_{\alpha}}\frac{d^2}{ds^2}\Bigg|_{s=0} \int_K \left( f(k\cdot \exp(sT_i^\alpha)a_t\cdot o)\right) Y_{\delta, j}(k)~dk
		\end{align*}
		The last expression, after the change of variable $k\mapsto k\cdot\exp(-sT_i^\alpha)$, becomes 
		\begin{align}
			\label{Thm1-eq-1}
			\sum_{i=1}^{m_{\alpha}} \int_K f(ka_t\cdot o)~ \frac{d^2}{ds^2}\Bigg|_{s=0} Y_{\delta, j}(k\cdot \exp(-sT_i^\alpha))~dk
		\end{align}
		But recalling the definition of $Y_{\delta, j}$ (See \eqref{Sph-Harm-def}), we note that 
		\begin{align*}
			\sum_{i=1}^{m_{\alpha}}\frac{d^2}{ds^2}\Bigg|_{s=0} Y_{\delta, j}(k\cdot \exp(-sT_i^\alpha)) &= \sum_{i=1}^{m_{\alpha}}\frac{d^2}{ds^2}\Bigg|_{s=0} (v_j,~ \delta (k)\delta(\exp(-sT_i^\alpha))v_1)_{L^2}\\
			& = \left(v_j,~ \delta (k)~ \sum_{i=1}^{m_{\alpha}}\frac{d^2}{ds^2}\Bigg|_{s=0}\delta(\exp(-sT_i^\alpha))v_1\right)_{L^2}\\
			& = (v_j, \delta (k)~\delta(\omega_{\alpha})v_1)_{L^2}, 
		\end{align*}
		which then from the above equations yields 
		\begin{align*}
			\int_K ((\Ad(a_t^{-1})\omega_{\alpha})f)(ka_t) Y_{\delta,j}(k)~dk= \int_K f(ka_t\cdot o)(v_j, \delta (k)~\delta(\omega_{\alpha})v_1)_{L^2}~dk\quad\quad(\alpha\in \{\gamma,2\gamma\}).
		\end{align*}
		Therefore, by using Theorem~\ref{omega-eigenvalue}, we obtain
		\begin{align*}
			\int_K ((\Ad(a_t^{-1})\omega_{\alpha})f)(ka_t) Y_{\delta,j}(k)~dk=u(t)\times \begin{cases}\frac{p(p+m_{2\gamma}-1)-q(q+m_{\gamma}+m_{2\gamma}-1)}{2(m_{\gamma}+4m_{2\gamma})}~~~&\text{if}~~\alpha=\gamma\\
				-\frac{2p(p+m_{2\gamma}-1)}{m_{\gamma}+4m_{2\gamma}}~~~~~&\text{if}~~\alpha=2\gamma. \end{cases}
		\end{align*}  
		Plugging the above into \eqref{sph-coeff-laplacef}, we then have 
		\begin{align*}
			&\int_{K} \Delta_{X}f(ka_t\cdot o)Y_{\delta,j}(k) ~dk
			= \frac12(m_\gamma+4m_{2\gamma})^{-1} \left(\frac{d^2}{d t^2} + (m_{\gamma} \coth t+ 2m_{\gamma} \coth (2t))\frac{d}{d t}\right)u(t) \\ &+\left(\frac{p(p+m_{2\gamma}-1)-q(q+m_{\gamma}+m_{2\gamma}-1)}{2(m_{\gamma}+4m_{2\gamma})}\sinh^{-2}(t)-\frac{2p(p+m_{2\gamma}-1)}{m_{\gamma}+4m_{2\gamma}}\sinh^{-2}(2t)\right)u(t)
		\end{align*}
		which, using \eqref{sph-coeff-Helm}, shows that $u(t)$ solves 
		\begin{align}
			\label{radial-equn-Helm}
			&\nonumber\left(\frac{d^2}{d t^2} + (m_{\gamma} \coth t+ 2m_{\gamma} \coth (2t))\frac{d}{d t}+  {\big(p(p+m_{2\gamma}-1)-q(q+m_{\gamma}+m_{2\gamma}-1)\big)}\sinh^{-2}(t)\right.\\ &\hspace*{2cm}\left.-{4p(p+m_{2\gamma}-1)}\sinh^{-2}(2t)+2(m_{\gamma}+4m_{2\gamma})(\lambda^2+\rho^2)\right)u(t)=0.
		\end{align}
		Now writing $l=i\lambda-\rho$, and recalling $2\rho=m_\gamma+2m_{2\gamma}$, by a straightforward calculation (see \cite[p.~543]{Vel} for a similar computation), we rewrite \eqref{radial-equn-Helm} as
		\begin{align}
			\label{radial-equn-Helm2}
			&	\frac1{(\cosh t)^{m_{2\gamma}}(\sinh t)^{m_{\gamma}+m_{2\gamma}}}~\frac{d}{dt} \left(	(\cosh t)^{m_{2\gamma}}(\sinh t)^{m_{\gamma}+m_{2\gamma}}\frac{du}{dt}\right)\nonumber \\
			&+\left(\frac{p(p+m_{2\gamma}-1)}{\cosh^2t}-\frac{q(q+m_{\gamma}+m_{2\gamma}-1)}{\sinh^2t}-l(l+m_{\gamma}+2m_{2\gamma})\right)u(t)=0
		\end{align}
		The above equation can be transformed into a standard hypergeometric differential equation through suitable substitutions. Specifically, set
		$$
		v(t) = \tanh^{-q}(t)\,\cosh^{-l}(t)\,u(t), 
		\quad\text{and then}\quad z = \tanh^2(t).
		$$
		After a straightforward (though somewhat lengthy) computation, equation \eqref{radial-equn-Helm2} takes the form
		\begin{equation}
			\label{hypergeometric-eq-under-consi}
			\frac{d^2 v}{dz^2} + \big[c - (a+b+1)z\big]\frac{dv}{dz} - ab\,v(z) = 0, 
			\quad\quad (z > C > 0),
		\end{equation}
		where
		\begin{align}
			\label{def:abc}
			\begin{cases}
				a = \tfrac12(q-l-p), \\[6pt]
				b = \tfrac12\big(q-l-p-m_{2\gamma}+1\big), \\[6pt]
				c = q + \tfrac12\big(m_{\gamma}+m_{2\gamma}+1\big).
			\end{cases}
		\end{align}
		{Since $z$ is bounded away from the origin, and $c-a-b=i\lambda\notin \mathbb{Z}$ in view of the hypothesis,  two linearly independent solutions of this hypergeometric equation are given by (See \cite[Ch.~5, Sec.~10.4]{Olver})}
		$$
		v_1(z) = {}_2F_1(a,b; a+b-c+1; 1-z),
		\quad 
		v_2(z) = (1-z)^{c-a-b}{}_2F_1(c-a, c-b; c-a-b+1; 1-z),
		$$
		where ${}_2F_1$ denotes the Gauss hypergeometric function. For more about hypergeometric functions, we refer the reader to \cite[Ch.~5, Sec.~9]{Olver}, see also \cite[Ch.~15]{OM}.

		Transforming back to the variable $t$, and recalling that $l=i\lambda-\rho$,  we conclude that for $t > R_0$ there exist constants $c_1, c_2$ such that
		$$
		u(t) = c_1 u_1(t) + c_2 u_2(t),
		$$
		where $u_1$ and $u_2$ are explicitly given by
		\begin{align}
			\label{exp:u1}
			u_1(t) &:= (\tanh t)^q(\cosh t)^{i\lambda-\rho} ~
			{}_2F_1\left(\tfrac{q-l+p}{2}, \tfrac{q-l-p-m_{2\gamma}+1}{2}; 1-i\lambda; 1-(\tanh t)^2\right),
		\end{align}
		and
		\begin{align}
			\label{exp:u2}
			u_2(t) &:= (\tanh t)^{q} (\cosh t)^{-i\lambda-\rho} \nonumber\\
			&\times{}_2F_1\left(\tfrac{p+l+q-m_\gamma+m_{2\gamma}+1}{2},  \tfrac{m_{\gamma}+2m_{2\gamma}+p+q+l}{2};  1+i\lambda;  1-(\tanh t)^2\right).
		\end{align}
		Here we have used the relation $c-a-b = i\lambda$, which follows from \eqref{def:abc}.

		Using the facts that $\tanh t \to 1$, $\cosh t \sim \tfrac{1}{2} e^{t}$ as $t \to \infty$, and the property of the hypergeometric function ${}_2F_1(a,b;c;0)=1$, we obtain 
		\begin{align*}
			u_1(t) &= e^{(i\lambda+\rho)t}\bigl(1 + o(1)\bigr), 
			& u_2(t) &= e^{(-i\lambda+\rho)t}\bigl(1 + o(1)\bigr),
			\quad (t \to \infty).
		\end{align*}
		
		Consequently, 
		\begin{align}
			\label{behv:u}
			u(t)=e^{-\rho t}\big(C_1e^{i\lambda t}+C_2e^{-i\lambda t}\big)~(1+o(1)) \quad\quad(t\rightarrow\infty).
		\end{align}
		For convenience, we introduce the notation
		$$g_{\lambda}(t)=C_1e^{i\lambda t}+C_2e^{-i\lambda t}\quad\quad (\lambda\in\mathbb{C}, t>R_0).$$
		From the above estimate \eqref{behv:u} it follows that for $t>R_0$, there exist $C>0$ such that 
		\begin{align}
			\label{u(t)-lowbdd-1}
			|u(t)|\geq C e^{-\rho t} |g_{\lambda}(t)|.
		\end{align}
		Therefore, recalling the definition of $u(t)$ from \eqref{def-u}, for any $1\leq p\leq\infty$, applying H\"older inequality (with standard modification for the $p=\infty$ case), we see that  
		\begin{align*}
			C e^{-\rho t} |g_{\lambda}(t)|\leq  \left| \int_{K} f(ka_t\cdot o)Y_{\delta,j}(k) ~dk\right|	\leq \left(\int_{K}|f(ka_t)|^p~dk\right)^{\frac1p}\left(\int_K|Y_{\delta,j}(k)|^{p'}~dk\right)^{\frac{1}{p'}}
		\end{align*}
		and subsequently, using the fact that $\|Y_{\delta,j}\|_{L^{p'}(K)}$ is always non-zero, we note that there exists a constant $C_3>0~(C_3\equiv C_3(\delta, j, p))$  such that 
		\begin{align}
			\label{main-estimate-ut}
			e^{-\rho t} |g_{\lambda}(t)|\leq C_3 \int_{K}|f(ka_t)|^p~dk\quad\quad(t>R_0). 
		\end{align}
		In what follows, we analyze the behavior of the function $|g_{\lambda}(t)|$ in different regimes depending on the value of $\lambda.$\\
		\underline{\textit{Case I}: $\Im(\lambda)\neq 0.$}\\
		We first consider the case where $\Im(\lambda)\neq0.$ First note that 
		$$g_{\lambda}(t)=C_1e^{i \Re(\la)t}e^{-\Im(\la)t}+C_2e^{-i \Re(\la)t}e^{\Im(\la)t}.$$ Now if either $C_1$ or $C_2$ is zero, the it is easy to see that
		\begin{align}
			\label{trivialcaseIMlaneq0}
			|g_{\lambda}(t)|=|C_2|e^{\Im(\la)t}~~\text{or}~~|C_1|e^{-\Im(\la)t},
		\end{align} respectively. Therefore, we now assume that both $C_1$, and $C_2$ are non-zero.

		We next first consider the case of $\Im(\lambda)>0$.  Then, $e^{\Im(\la)t}$ being the dominating term in $g_{\lambda}(t)$ for large $t$, we observe that 
		$$|g_{\lambda}(t)|\sim |C_2|e^{\Im(\la)t}\quad\quad (t\rightarrow\infty).$$ 
		Similarly, when $\Im(\la)<0$, the dominating factor becomes $e^{-\Im(\la)t}$, and subsequently, 
		$$|g_{\lambda}(t)|\sim |C_1|e^{-\Im(\la)t}\quad\quad (t\rightarrow\infty).$$ 
		Therefore, using above observations along with \eqref{trivialcaseIMlaneq0}, we conclude that for $\Im(\la)\neq 0$, there exists $C^{\pm}>0$ (independent of $t$), such that 
		\begin{align}
			\label{asymp-g-la-t}
			|g_{\lambda}(t)|\sim C^{\pm}e^{\pm\Im(\la)t}\quad\quad (t\rightarrow\infty).
		\end{align}
		Consequently,  using the above asymptotics, there exists $R_1>R_0$ (sufficiently large), and a constant $C>0$ such that for all $t>R_1$ we have  
		\begin{align}
			|g_{\la}(t)|\geq C e^{-| \Im(\la)|t}. 
		\end{align}
		Therefore, in view of \eqref{main-estimate-ut}, for any $1\leq p\leq\infty$, {with standard modification for $p=\infty$}, we have 
		\begin{align}
			e^{p(-|\Im(\la)|-\rho)t}\leq C_3 \int_{K}|f(ka_t)|^p~dk\quad\quad(t>R_1). 
		\end{align}
		\underline{\textit{Case II}: $\Im(\lambda)=0.$}\\
		
		A straightforward calculation involving basic complex analysis gives  
		\[
		|g_{\la}(t)|^2=\Big|C_1 e^{i\lambda t} + C_2 e^{-i\lambda t} \Big|^2 = |C_1|^2 + |C_2|^2 + 2 |C_1| |C_2| \cos(2\lambda t + \theta),
		\]
		where \(\theta = \arg(C_1 \overline{C_2})\). 
		This, in view of \eqref{main-estimate-ut}, implies  
		\[
		Ce^{-2\rho t}\Big(|C_1|^2 + |C_2|^2 + 2 |C_1| |C_2| \cos(2\lambda t + \theta)\Big)\leq \left(\int_K|f(ka_t)|^p\right)^\frac2p
		\]
		completing the proof of the proposition.
	\end{proof}
	\begin{rem}
		\label{rem:exclusion-dis}
		The technical restriction $\lambda\notin i\mathbb{Z}$ on the spectral parameter $\lambda$ arises directly from the general 
		theory of the hypergeometric differential equation. 
		As seen in the above proof, this exclusion is necessary to ensure the linear independence 
		of the two fundamental solutions $v_1$ and $v_2$ of the radial differential equation~\eqref{hypergeometric-eq-under-consi} 
		near one of its regular singular points, namely $z=1$, which is crucial in our analysis.

		More precisely, the local behavior near \(z=1\) depends on the quantity \(c-a-b\). When \(c-a-b \notin \mathbb{Z}\), the two solutions are linearly independent and exhibit standard power-type behavior. In contrast, when \(c-a-b \in \mathbb{Z}\), a resonance occurs: the second solution develops a logarithmic term, leading to a loss of analytic independence and more complicated asymptotics. Such cases are therefore excluded, as they obstruct the derivation of the required \(L^{p}\)-estimates.
		
		For further details on the hypergeometric equation near regular singular points, see \cite[Chap.~15]{OM} and \cite{Olver}.
		\end{rem}
	We now prove Theorem~\ref{rellich-asymptotic-rss}.
	\begin{proof}[Proof of Theorem \ref{rellich-asymptotic-rss}:]
		We divide the argument into two cases, according to whether the spectral parameter 	$\lambda$ is real or has a nonzero imaginary part.
		\medskip
		\noindent
		We first consider the case $\Im(\la)\neq 0.$ Let $1\leq p<2$.   
		Using the polar decomposition of the measure \eqref{polar-int}, and the estimate \eqref{est-la-non-0} of Proposition \ref{main-prop},  there exists $R_{1}>R_0$ such that for any $R>R_1$, we get
		\begin{align*}
			\int_{R<d(o,~x)<2R}|f(x)|^p~dx= \int_R^{2R}J(t) \int_{K}|f(ka_t)|^p~dk~dt\geq C \int_{R}^{2R} e^{p(-|\Im(\la)|-\rho)t} e^{2\rho t} ~dt.
		\end{align*}
		Recalling $\gamma_p=\tfrac{2}{p}-1$, it is easy to observe that 
		$$p(-|\Im(\la)|-\rho)+2\rho= p(\gamma_p\rho-|\Im(\la)|).$$
		But then  when $|\Im(\la)|\neq \gamma_p\rho$, the integral on the right can be estimated as 
		\begin{align*}
			I=&	\int_{R}^{2R} e^{p(-|\Im(\la)|-\rho)t} e^{2\rho t} ~dt\\
			&=\int_{R}^{2R} e^{(-|\Im(\la)|+\gamma_p\rho)pt}  ~dt\\
			&= \frac{1}{(\gamma_p\rho-|\Im(\la)|)p}\left(e^{2p(\gamma_p\rho-|\Im(\la)|)R}-e^{p(\gamma_p\rho-|\Im(\la)|)R}\right)\\
			&= \frac{1}{(\gamma_p\rho-|\Im(\la)|)p}e^{p(\gamma_p\rho-|\Im(\la)|)R} \left(e^{p(\gamma_p\rho-|\Im(\la)|)R}-1\right)\\
			&= e^{p(\gamma_p\rho-|\Im(\la)|)R} L(\la,R,p)
		\end{align*}
		where notice that  
		\begin{align*}
			L(\la,R,p):= \begin{cases}
				\frac{\left(e^{p(\gamma_p\rho-|\Im(\la)|)R}-1\right)}{(\gamma_p\rho-|\Im(\la)|)p}&\text{when}~~|\Im(\la)|< \gamma_p\rho \\
				\frac{\left(1-e^{p(\gamma_p\rho-|\Im(\la)|)R}\right)}{(|\Im(\la)|-\gamma_p\rho)p}& \text{when}~~|\Im(\la)|> \gamma_p\rho
			\end{cases}
		\end{align*}
		Since the factor $e^{p(\gamma_p \rho - |\Im(\lambda)|)R}$ grows rapidly as $R \to \infty$ when 
		$|\Im(\lambda)| < \gamma_p \rho$, while it decays rapidly when $|\Im(\lambda)| > \gamma_p \rho$, 
		we can choose $R_1 > R_0$ sufficiently large so that for all $R > R_1$,
		\[
		L(\lambda, R, p) \geq C(p, \lambda).
		\]
		Consequently,
		\[
		I \geq C(p, \lambda) \, e^{p(\gamma_p \rho - |\Im(\lambda)|)R}.
		\]
		Hence, for any $R > R_1 $, we obtain
		\[
		\int_{R < d(o,\,x) < 2R} |f(x)|^p \, dx 
		\;\geq\;
		C \, e^{p(\gamma_p \rho - |\Im(\lambda)|)R},
		\]
		where the constant $C$ may depend on $\lambda$, $p$, and $\rho$, but not on $R$.
		This completes the proof in the present case.

		Next, we consider the critical regime $|\Im(\lambda)| = \gamma_p \rho$. 
		Here the exponential factor vanishes, leaving  
		\[
		I = \int_R^{2R} e^{(-|\Im(\lambda)| + \gamma_p \rho)pt}\, dt = \int_R^{2R} 1\, dt = R.
		\]
		Hence,  
		\[
		\int_{R < d(o,\,x) < 2R} |f(x)|^p \, dx \;\geq\; C R.
		\] 
		Next, we  assume that $\Im(\la)=0.$  By part (2) of Proposition \eqref{main-prop}( with $p=2$), for any $R>R_0$ we have
		\begin{align*}
			C\int_R^{2R} \big(A^2+B^2+2AB\cos(2\lambda t+\theta)\big)~dt&\leq \int_R^{2R}e^{2\rho t}\int_K|f(ka_t)|^2~dk~dt\\&\leq\int_{R<d(o,~x)<2R} |f(x)|^2~dx.
		\end{align*} 
		The integral on the left can be evaluated explicitly:
		\begin{align*}
			\int_R^{2R} \big(A^2+B^2+2AB\cos(2\lambda t+\theta)\big)~dt= (A^2+B^2) R - \frac{AB}{\lambda} \left( \sin(2\lambda R + \theta) - \sin(4\lambda R + \theta) \right).
		\end{align*} 
		Choosing \(R_1 > R_0\) sufficiently large, there exists \(C > 0\) such that  for all $R>R_1$ we have 
		\[
		(A^2+B^2) R - \frac{AB}{\lambda} \left( \sin(2\lambda R + \theta) - \sin(4\lambda R + \theta) \right) > C R.
		\]
		Thus,  for all $R>R_1$, we arrive at 
		\[
		\int_{R < |x| < 2R} |f(x)|^2 dx \geq C R.
		\]
		This completes the proof of part (1), and hence the theorem.
	\end{proof}
	
	We now employ estimates in Theorem~\ref{rellich-asymptotic-rss} to establish Theorem~\ref{rel-main-thm}. To this end, we shall make use of a Pizzetti-type mean value property for analytic functions defined on open subsets of $X$ (see \cite[(4.1)]{Q}; see also \cite[(7.5)]{BZ}).
	
	\begin{lem}
		\label{mvt}
		Let $U\subset X$ be open and connected, and $f$ be an analytic function on $U.$ Then for sufficiently small $r>0$, $f$ satisfies
		\begin{equation}
			\frac{1}{|S_r(x)|}\int_{S_r(x)}f(y)~d\sigma(y)=\Gamma(\tfrac{n}{2})\sum_{m=0}^\infty \left(\frac{\sinh \kappa r}{2\kappa}\right)^{2m} L_mf(x),
		\end{equation} 
		where $d\sigma$ is the measure on the sphere $S_t(x)$, and the operator $L_m$ is defined by $$L_m=\frac{\Delta_{X}(\Delta_{X}-(2n+4\beta+4)\kappa^2)\cdots (\Delta_{X}-(2n+4\beta+4m-4)\kappa^2}{m!\Gamma(m+\tfrac{n}{2})}.$$
		Here $\beta$, and $\kappa$ are constants related to the structure and curvature of $X$, respectively.
	\end{lem}
	\begin{proof}[Proof of Theorem~\ref{rel-main-thm}]
		Let $f$ satisfy \eqref{helmrss} in $\Omega$, and suppose that $f \in L^{p}(\Omega)$. We argue by contradiction and assume that $f \neq 0$.
		
		First, consider the case $p = 2$. In this case, $\gamma_{p} = 0$, and consequently, by the part~(1) of Theorem~\ref{rellich-asymptotic-rss}, $f$ satisfies \eqref{rel-l2-est-rss}, which contradicts the assumption that $f \in L^{2}(\Omega)$.
		
		Next, assume that $1 \leq p < 2$. First, let $\Im(\lambda) \neq 0$. Under this hypothesis on $\lambda$, we have $
		\gamma_{p}\rho - |\Im(\lambda)| > 0.$
		In view of part~(2) of Theorem~\ref{rellich-asymptotic-rss}, this implies that $f$ satisfies \eqref{rel-lp-est-rss}, which again contradicts the assumption that $f \in L^{p}(\Omega)$.
		
		Finally, we consider the remaining case $\Im(\lambda) = 0$. To obtain a contradiction, we employ the $L^{2}$ estimate \eqref{rel-l2-est-rss} together with the mean-value type result in Lemma~\ref{mvt}. 
		To begin with, note that since $-\Delta_X$ is elliptic, the function $f$ is real-analytic on any open, connected subset $U \subset \Omega$. In particular, $f \not\equiv 0$ on $U$, and hence 
		\[
		\|f\|_{L^\infty(\Omega_1)} > 0, 
		\qquad 
		\Omega_1 := \{x \in X : d(o, x) > R_0 + 1\} \subset \Omega.
		\]
		We next show that $\|f\|_{L^\infty(\Omega_1)} < \infty$.
		
		Let $x \in \Omega_1$. Since $f$ is an eigenfunction of $\Delta_X$ with eigenvalue $-(\lambda^2 + \rho^2)$ in $\Omega$, there exists $r_1 > 0$ sufficiently small such that, for all $0 < t < r_1$, Lemma~\ref{mvt} yields the mean value identity
		\begin{align}
			\label{mvt-f}
			\frac{1}{|S_t(x)|} \int_{S_t(x)} f(y)\, d\sigma(y)
			= f(x)\, \Gamma\!\left(\tfrac{n}{2}\right)
			\sum_{m=0}^\infty 
			\left(\frac{\sinh(\kappa t)}{2\kappa}\right)^{2m}
			\frac{\Psi_m(\lambda, n, \beta, \rho)}{m!\, \Gamma(m + \tfrac{n}{2})},
		\end{align}
		where
		\[
		\Psi_m(\lambda, n, \beta, \rho)
		= (-1)^m (\lambda^2 + \rho^2)
		\big[(\lambda^2 + \rho^2) + (2n + 4\beta + 4)\kappa^2\big]
		\cdots
		\big[(\lambda^2 + \rho^2) + (2n + 4\beta + 4m - 4)\kappa^2\big].
		\]
		
		Define
		\[
		\Psi(t)
		:= \Gamma\!\left(\tfrac{n}{2}\right)
		\sum_{m=0}^\infty 
		\left(\frac{\sinh(\kappa t)}{2\kappa}\right)^{2m}
		\frac{\Psi_m(\lambda, n, \beta, \rho)}{m!\, \Gamma(m + \tfrac{n}{2})}.
		\]
		Then $\Psi$ is continuous and satisfies $\Psi(0) = 1$. 
		Since both $\Psi$ and the Jacobian $J$ of the polar decomposition   $J$ are continuous and $J(t) > 0$ for $t>0$ (see~\eqref{def-J}), one can choose $0 < r < r_1$ small enough so that
		\[
		D(r):= \int_0^r \Psi(t)\, J(t)\, dt \neq 0.
		\]
		(The dependence of $\Psi$ and $D(r)$ on the parameters appearing in $\Psi_m$ is immaterial and will be suppressed in the notation.)
		
		Now, recalling the fact that 
		\[\frac{1}{|S_t(x)|} \int_{S_t(x)} f(y)\, d\sigma(y)=\int_{K}f(xka_t)~dk,\]
		from \eqref{mvt-f}, we see that
		\begin{align*}
			\frac{1}{|B_r(x)|} \int_{B_r(x)} f(y)\, dy=\frac{1}{|B_r(x)|}\int_{B_r(o)}f(xy)~dy=\int_0^r\int_{K} f(xka_t)J(t)~dk~dt=f(x)\, D(r),
		\end{align*}
		where in the second equality we have used the polar decomposition of the measure (See \eqref{polar-int}).
		Hence, by Hölder’s inequality,
		\[
		|f(x)|
		\leq \frac{1}{|D(r)|\, |B_r(o)|} \int_{B_r(x)} |f(y)|\, dy
		\leq \frac{|B_r(o)|^{\tfrac{1}{p'}-1}}{|D(r)|} 
		\left( \int_{\Omega} |f(y)|^p\, dy \right)^{1/p},
		\]
		and therefore $\|f\|_{L^\infty(\Omega_1)} < \infty$.
		
		Combining this bound with the $L^2$-estimate~\eqref{rel-l2-est-rss} yields, for all $R > \max(R_1, R_0 + 2)$,
		\begin{equation}
			\label{L2-Lp-est}
			C R
			\leq \int_{R < d(o, x) < 2R} |f(x)|^2\, dx
			\leq
			\|f\|_{L^\infty(\Omega_1)}^{2 - p}
			\int_{R < d(o, x) < 2R} |f(x)|^p\, dx,
		\end{equation}
		contradicting the assumption that $f \in L^p(\Omega)$. 
		This establishes the claim when $\Im(\lambda) = 0$, and hence completes the proof of the theorem.
	\end{proof}
	\begin{rem}
		\label{rem:sharpness-lp-est}
		To examine the sharpness of Theorem~\ref{rel-main-thm}, we construct nontrivial eigenfunctions on 
		$\Omega$ with eigenvalue $-(\lambda^2+\rho^2)$ satisfying the $L^p$-integrability condition~\eqref{rel-X-Lp-est} 
		for pairs $(p,\lambda)$ with $\lambda\notin i\mathbb{Z}$ not covered by Theorem~\ref{rel-lp-est-rss}.
		
		\smallskip
		\noindent
		\textbf{(i) The case $p>2$, interior of the $L^p(X)$ spectrum.}
		It is well known (see~\eqref{phi-la-lp-est}) that the elementary spherical functions 
		$\varphi_{\lambda}\in L^p(X)$ whenever 
		\[
		p>2 \quad\text{and}\quad |\Im(\lambda)|<(1-\tfrac{2}{p})\rho.
		\]
		Hence, for such parameters we may simply take $f=\varphi_{\lambda}$ as an example.
		
		\smallskip
		\noindent
		\textbf{(ii) The complementary region.}
		For the remaining ranges of $(p,\lambda)$, namely
		\[
		1\le p<2,\ |\Im(\lambda)|>\gamma_p\rho, 
		\quad\text{or}\quad 
		p>2,\ |\Im(\lambda)|\ge (1-\tfrac{2}{p})\rho,
		\]
		we use the function $\Phi_{\lambda}$ appearing in the Harish--Chandra expansion of $\varphi_{\lambda}$ 
		(see~\eqref{HC-expan}). From the asymptotic estimate~\eqref{Phi-asymp}, for large $t$ and any 
		$1\le p<\infty$,
		\begin{equation}
			\label{Phi-K-Lp-mean}
			\int_{K}|\Phi_{\lambda}(ka_t)|^p\,dk \asymp e^{p(-\Im(\lambda)-\rho)t}.
		\end{equation}
		Consequently, recalling $\gamma_p = \tfrac{2}{p}-1$,
		\begin{align*}
			\int_{\Omega}|\Phi_{\lambda}(x)|^p\,dx
			&\asymp 
			\int_{R_0}^{\infty} e^{p(-\Im(\lambda)-\rho)t} e^{2\rho t}\,dt
			= \int_{R_0}^{\infty} e^{-p(\Im(\lambda)-\gamma_p\rho)t}\,dt
			= \int_{R_0}^{\infty} e^{-p\nu(\lambda,p)t}\,dt,
		\end{align*}
		where $\nu(\lambda,p):=\Im(\lambda)-\gamma_p\rho$.
		
		Hence $\Phi_{\lambda}\in L^p(\Omega)$ whenever $\nu(\lambda,p)>0$.  
		For $p>2$ (so that $\gamma_p<0$), this condition is automatically satisfied when $\Im(\lambda)>0$, 
		and we take $f=\Phi_{\lambda}$. 
		If $\Im(\lambda)<0$, we instead choose $f=\Phi_{-\lambda}$.  
		Similarly, for $1\le p<2$, $L^p$–integrability over $\Omega$ requires $|\Im(\lambda)|>\gamma_p\rho$, 
		and we again select $f=\Phi_{\pm\lambda}$ according to the sign of $\Im(\lambda)$.
		
		\smallskip
		\noindent
		These constructions show that in all complementary spectral regions, one can find 
		nontrivial $L^p$–eigenfunctions of $-\Delta_X$ on $\Omega$, thereby confirming the sharpness of 
		Theorem~\ref{rel-main-thm}.
	\end{rem}
	
	We now turn to the proof of Theorem~\ref{rel-rss-weak-lp}. The argument combines the \(L^{p}\)-\(K\)-mean estimates from Proposition~\ref{main-prop} with a standard measure-theoretic characterization of weak \(L^{p}\)-spaces, which allows us to pass from averaged estimates to weak \(L^{p}\)-bounds.
	
	\begin{lem}
		Let \((M,\mu)\) be a measure space and \(0<p<\infty\). Fix \(0<r<p\), and define
		\begin{equation}
			\label{def:new-lp}
			|||f|||_{L^{p,\infty}(M)}
			:=\sup_{0<\mu(E)<\infty}\mu(E)^{-\frac1r+\frac1p}
			\left(\int_E|f(x)|^r\,d\mu(x)\right)^{\frac1r},
		\end{equation}
		where the supremum is taken over all measurable sets \(E \subset M\) of finite measure. Then, for \(f \in L^{p,\infty}(M)\),
		\begin{equation}
			\label{norm-eqv-weak-lp}
			\|f\|_{L^{p,\infty}(M)}
			\leq |||f|||_{L^{p,\infty}(M)}
			\leq \left(\frac{p}{p-r}\right)^{\frac1r}\|f\|_{L^{p,\infty}(M)}.
		\end{equation}
	\end{lem}
	
	We refer the reader to \cite[1.1.12, pp.14]{Grafakos} for this result. We also remark in passing that, for \(p>1\), the quantity \(|||\cdot|||\) defines a norm equivalent to the usual weak \(L^{p}\)-(quasi-) norm. This equivalence will be used repeatedly in what follows to convert integral estimates into weak \(L^{p}\)-bounds.
	
	\begin{proof}[Proof of Theorem~\ref{rel-rss-weak-lp}]
		Fix \(1<p<2\) and \(\lambda \in \mathbb{C}\setminus i\mathbb{Z}\). Assume that \(f\not\equiv 0\) solves~\eqref{helmrss} in \(\Omega\).
		
		We first prove part~(1). Suppose \(\lambda=\alpha\pm i\gamma_p\rho\) with \(\alpha\in\mathbb{R}\), and \(f \in L^q(\Omega)\) for some \(1 \leq q \leq p\). By Proposition~\ref{main-prop}, there exists \(R_1>R_0\) such that
		\[
		\int_K |f(ka_t)|^q\,dk \gtrsim e^{-(\rho+\gamma_p\rho)qt}, \qquad t>R_1.
		\]
		Hence,
		\begin{align*}
			\|f\|_{L^q(\Omega)}^q
			&\geq \int_{R_1}^\infty \int_K |f(ka_t)|^q e^{2\rho t}\,dk\,dt \\
			&\gtrsim \int_{R_1}^\infty e^{(2\rho-\rho q-\gamma_p\rho q)t}\,dt.
		\end{align*}
		Since \(2\rho-\rho q-\gamma_p\rho q = 2\rho(1-q/p)\geq 0\), the last integral diverges, a contradiction. Thus \(f=0\).
		
		Next, we prove part~(2). Assume \(f \in L^{q,\infty}(\Omega)\) with \(1<q<p\). For \(R>0\), set \(E_R=(R,2R)\times K\). Then, in view of the integraion in polar co-ordinates \eqref{polar-int}, we record that 
		\[
		|E_R|\asymp\int_R^{2R} e^{2\rho t}\,dt =\frac{1}{2\rho}\left(e^{4\rho R}-e^{2\rho R}\right).
		\]
		Using~\eqref{def:new-lp} (with \(r=1\)), we obtain
		\[
		|||f|||_{L^{q,\infty}(\Omega)}
		\gtrsim |E_R|^{-1+\frac1q}\int_{E_R} |f(ka_t)| e^{2\rho t}\,dk\,dt.
		\]
			But, then using the estimate \eqref{rel-lp-est-rss} of the Proposition~\ref{main-prop}, there exist $R_1>R_0$ such that for all $R>R_1$, we have $$\int_{E_R}|f(ka_t)|e^{2\rho t}~dkdt\gtrsim \int_R^{2R}e^{-\rho t-|\Im(\lambda)|t}e^{2\rho t}dt=\frac{ e^{2R(\rho-|\Im(\lambda)|)}-e^{R(\rho-|\Im(\lambda)|)}}{\rho-|\Im(\lambda)|}.$$
			Therefore, in view of \eqref{norm-eqv-weak-lp}, for any $R>R_1$, we obtain \begin{equation} \label{key-est-thm-weak-lp} \|f\|_{L^{q,\infty}(\Omega)}\gtrsim \left(e^{4\rho R}-e^{2\rho R}\right)^{\frac1q-1}\left(e^{2R(\rho-|\Im(\lambda)|)}-e^{R(\rho-|\Im(\lambda)|)}\right). \end{equation} Now, if $\lambda=\alpha\pm i\gamma_p\rho$, and $1<q<p$, in view of the fact that $\rho-|\Im(\lambda)|=2\rho(1-1/p)$, using a simple calculation we see from \eqref{key-est-thm-weak-lp} that
		\[
		\|f\|_{L^{q,\infty}(\Omega)} \gtrsim e^{4\rho R(\frac1q-\frac1p)} \to \infty
		\quad \text{as } R\to\infty,
		\]
		since \(q<p\). This contradicts the assumption \(f \in L^{q,\infty}(\Omega)\), and thus \(f=0\).
		
		Finally, we prove part~(3). Let \(\lambda=\alpha\pm i\gamma_q\rho\) with \(p<q< 2\), and assume \(f \in L^{p,\infty}(\Omega)\). Repeating the above argument, we obtain for large \(R\)
		\[
		\|f\|_{L^{p,\infty}(\Omega)} \gtrsim e^{4\rho R(\frac1p-\frac1q)} \to \infty,
		\]
		which is again a contradiction. This completes the proof.
	\end{proof} 
\begin{rem}
	\label{weak-lp-sharpness}
	The result is sharp in the following sense. There exist nontrivial eigenfunctions in \(\Omega\) belonging to \(L^{q,\infty}(\Omega)\) for \(q \geq p\) (resp. \(L^{p,\infty}(\Omega)\) for \(q \leq p\)) when the spectral parameter is of the form \(\lambda = \alpha \pm i \gamma_p \rho\) (resp. \(\lambda = \alpha \pm i \gamma_q \rho\)).
	
	To see this, we employ the functions \(\Phi_{\pm\lambda}\). We begin with a simple observation. Let \(A>0\), and consider the \(K\)-invariant function on \(\Omega\) defined by \(f(a_t)=e^{-At}\) for \(t>R_0\). By \(K\)-invariance,
	\[
	\{x \in \Omega : |f(x)| > s\}
	= \{t > R_0 : e^{-At} > s\}
	= \left(R_0, \tfrac{1}{A}\log \tfrac{1}{s}\right).
	\]
	Thus, recalling the definition of the  distribution function \(d_f(s)\) of $f$, we have \(d_f(s)=0\) for \(s \geq 1\), and for \(0<s<1\),
	\begin{align*}
		d_f(s)
		\asymp \int_{R_0}^{\frac{1}{A}\log \frac{1}{s}} e^{2\rho t}\,dt
		\lesssim s^{-\frac{2\rho}{A}}.
	\end{align*}
	Consequently, for any \(q>0\),
	\begin{equation}
		\label{weak-sharpness-gen-est}
		\|f\|_{L^{q,\infty}(\Omega)}
		= \sup_{0<s<1} s\, d_f(s)^{1/q}
		\lesssim \sup_{0<s<1} s^{1-\frac{2\rho}{qA}} < \infty
		\quad \text{whenever } qA \geq 2\rho.
	\end{equation}
	
	Now, using the estimate~\eqref{Phi-la-est}, we have for large $t$
	\[
	|\Phi_{\alpha + i\gamma_p \rho}(a_t)| \asymp e^{-\frac{2\rho}{p}t}.
	\]
	Applying~\eqref{weak-sharpness-gen-est}, it follows that \(\Phi_{\alpha + i\gamma_p \rho} \in L^{q,\infty}(\Omega)\) whenever \(q \geq p\). Similarly, \(\Phi_{\alpha + i\gamma_q \rho} \in L^{p,\infty}(\Omega)\) whenever \(q \leq p\). This establishes the sharpness.
\end{rem}
	\medskip

	Next, we turn our attention to the proof of Theorem~\ref{rel-hardy-norm}. This result follows immediately from Proposition~\ref{main-prop}. We begin by introducing the function $\psi_\lambda$  appearing in the statement of Theorem~\ref{rel-hardy-norm}. Its definition is naturally motivated by the lower bounds for the $L^p$–$K$–means of eigenfunctions established in Proposition~\ref{main-prop}. For $\lambda \in \mathbb{C}$, define the function $\psi_\lambda$ on 
	$\overline{A^+}=\{a_t : t \ge 0\}$ by
	\begin{equation}
		\label{def-psi}
		\psi_\lambda(a_t) = e^{(-|\Im(\lambda)| - \rho)t}, \qquad t \ge 0,
	\end{equation}
	and extend it to $G$ as a $K$–bi-invariant function. 
	The function $\psi_\lambda$ serves as a convenient comparison profile, 
	reflecting the principal exponential decay rate of eigenfunctions associated 
	with the spectral parameter $\lambda$. 
	Furthermore, in view of the estimate~\eqref{Phi-la} for $\Phi_{\lambda}$, we have the following 
	asymptotic equivalence for large~$t$:
	\begin{equation}
		\label{psi-asymp}
		\psi_{\lambda}(a_t) \asymp
		\begin{cases}
			|\Phi_{-i\,\Im(\lambda)}(a_t)|, & \text{if } \Im(\lambda) < 0, \\[4pt]
			|\Phi_{\,i\,\Im(\lambda)}(a_t)|, & \text{if } \Im(\lambda) > 0,
		\end{cases}
		\qquad \text{as } t \to \infty.
	\end{equation}
	
	\begin{proof}[Proof of Theorem~\ref{rel-hardy-norm}]
		Let $\lambda \in \mathbb{C}\setminus i\mathbb{Z}$, $\epsilon > 0$, $1 \leq p \leq \infty$, and let $f \in C^2(\Omega)$ be a solution of \eqref{helmrss} in $\Omega$ satisfying the growth condition \eqref{rel-hardy-norm-est}, i.e.,
		\begin{equation}
			\label{rel-hardy-norm-est-1}
			\sup_{t > R_0} 
			\frac{t^{\epsilon}}{\psi_{\lambda}(a_t)} 
			\left( \int_K |f(k a_t)|^p \, dk \right)^{\frac{1}{p}} < \infty.
		\end{equation}
		We argue by contradiction and assume that $f \not\equiv 0$. In what follows, we consider the case $p < \infty$; the case $p = \infty$ can be handled by a straightforward modification of the same argument.
		
		\medskip
		\noindent
		\textbf{Case 1.} $\boldsymbol{\Im(\lambda) \neq 0.}$  
		By part (1) of Proposition~\ref{main-prop}, there exists $R_1>R_0$ such that for all $t > R_1$, we have
		\[
		\int_K |f(k a_t)|^p \, dk \geq  C \, e^{p(-|\Im(\lambda)| - \rho)t}.
		\]
		This yields
		\[
		C t^{\epsilon} 
		\leq \frac{t^{\epsilon}}{e^{(-|\Im(\lambda)| - \rho)t}} 
		\left( \int_K |f(k a_t)|^p \, dk \right)^{\!\frac{1}{p}}
		= \frac{t^{\epsilon}}{\psi_{\lambda}(a_t)} 
		\left( \int_K |f(k a_t)|^p \, dk \right)^{\!\frac{1}{p}}.
		\]
		Consequently,
		\[
		\sup_{t > R_0} 
		\frac{t^{\epsilon}}{\psi_{\lambda}(a_t)} 
		\left( \int_K |f(k a_t)|^p \, dk \right)^{\!\frac{1}{p}} 
		= \infty,
		\]
		which contradicts \eqref{rel-hardy-norm-est-1}. Hence, we must have $f \equiv 0$ in $\Omega$.
		
		\medskip
		\noindent
		\textbf{Case 2.} $\boldsymbol{\Im(\lambda) = 0.}$  
		
		By part (2) of Proposition~\ref{main-prop}, for all $t>R_0$, we have
		\[
		\left( \int_K |f(k a_t)|^p \, dk \right)^{\!\frac{1}{p}}
		\geq C e^{-\rho t} U(\lambda, t),
		\]
		where
		\[
		U(\lambda, t) = \big(A^2 + B^2 + 2AB \cos(2\lambda t + \theta)\big)^{\frac{1}{2}}.
		\]
		Since $\psi_\lambda(a_t) = e^{-\rho t}$ whenever $\Im(\lambda) = 0$ (see \eqref{def-psi}), it follows that
		\[
		\frac{t^{\epsilon}}{\psi_{\lambda}(a_t)} 
		\left( \int_K |f(k a_t)|^p \, dk \right)^{\!\frac{1}{p}}
		\geq C t^{\epsilon} U(\lambda, t).
		\]
		However, it is easy to verify that
		\[
		\sup_{t > R_0} \big( t^{\epsilon} U(\lambda, t) \big) = \infty,
		\]
		and thus,
		\[
		\sup_{t > R_0}
		\frac{t^{\epsilon}}{\psi_{\lambda}(a_t)} 
		\left( \int_K |f(k a_t)|^p \, dk \right)^{\!\frac{1}{p}}
		= \infty,
		\]
		which contradicts \eqref{rel-hardy-norm-est-1}. Therefore, we conclude that $f \equiv 0$ in $\Omega$.
		\medskip
		Finally, we show that the condition $\epsilon>0$ in the hypothesis of the theorem is sharp. 
		To this end, we exhibit examples of eigenfunctions on $\Omega$ satisfying~\eqref{rel-hardy-norm-est-1} 
		with $\epsilon=0$.
		
		\smallskip
		\noindent
		First, let $\Im(\lambda)>0$. In this case, $\psi_{\lambda}(a_t)=e^{(-\Im(\lambda)-\rho)t}$, 
		and we choose $f=\Phi_{\lambda}$. 
		From the estimate~\eqref{Phi-K-Lp-mean}, it follows that
		\[
		\sup_{t > R_0} 
		\frac{1}{\psi_{\lambda}(a_t)} 
		\left( \int_K |\Phi_{\lambda}(k a_t)|^p \, dk \right)^{\!1/p} < \infty.
		\]
		Similarly, when $\Im(\lambda)<0$, taking $f=\Phi_{-\lambda}$ yields the same bound.
		
		\smallskip
		\noindent
		Finally, when $\Im(\lambda)=0$, we have $\psi_{\lambda}(a_t)=e^{-\rho t}$. 
		By Theorem~\ref{poisson-charac}, one may take $f=\mathcal{P}_{\lambda}F$ for any 
		$F\in L^p(K/M)$, which also satisfies~\eqref{rel-hardy-norm-est-1} with $\epsilon=0$.
		
		\smallskip
		\noindent
		Hence, from all the cases above, it follows that the assumption $\epsilon>0$ cannot be relaxed, 
		and the condition in the theorem is therefore sharp.
	\end{proof}
	\section{Concluding remarks}
	To close, we record a few observations and indicate some open questions that, in our opinion, 
	deserve further attention.
	
	We begin with the following observation. 
	As we have already seen, in the exterior domain $\Omega$, the Laplace--Beltrami operator 
	admits $L^p(\Omega)$ eigenfunctions with eigenvalue $-(\lambda^2+\rho^2)$ 
	for $1 \leq p \leq 2$ whenever $|\Im(\lambda)|>\gamma_p\rho$, 
	and for all $p>2$ with any $\lambda \in \mathbb{C}\setminus i\mathbb{Z}$, 
	owing to the presence of the spherical solutions $\Phi_{\pm\lambda}$. 
	Thus, as demonstrated throughout this article, a Rellich-type $L^p$ uniqueness result 
	is not possible for $p>2$. 
	However, in view of the $L^{2}$–estimate 
	\eqref{rel-l2-est-rss}, using H\"older's inequality, one can deduce a Rellich-type uniqueness statement 
	expressed in terms of appropriate weighted $L^{p}$–norms.
	
	\smallskip
	\noindent
	Let $\Omega$ be the exterior domain defined in~\eqref{def-Omega}, and suppose that 
	$f\in C^2(\Omega)$ is a nontrivial solution of the Helmholtz equation~\eqref{helmrss} 
	in $\Omega$ for $\lambda\in \mathbb{C}\setminus i\mathbb{Z}$. Then we have the following statement.
	\begin{itemize}
		\item 
		Let $1\leq p<2$, and let $\lambda=\alpha\pm i\gamma_p\rho$ with $\alpha\in\mathbb{R}$. 
		If 
		\[
		\int_{\Omega}|f(x)|^{p'}e^{2\gamma_p p'\rho\,d(o,x)}\,dx<\infty,
		\]
		then necessarily $f=0$ in $\Omega$.
	\end{itemize}
	
	\smallskip
	\noindent
	
	We conclude with a few natural questions that arise from the analysis carried out in this article.
	
	\begin{enumerate}[(1)]
		\item 
		As explained earlier, our method—based on the detailed analysis of the 
		corresponding hypergeometric differential equation—does not extend to the 
		spectral parameters $\lambda \in i\mathbb{Z}$. 
		It would be interesting to investigate what happens at these values, 
		in particular, to understand the precise behavior of the eigenfunctions 
		corresponding to such parameters and whether analogous Rellich-type estimates 
		still hold or fail.
		
		\item 
		A crucial step in our proof is the use of the spherical harmonic expansion, 
		where the  $K$–part (or rotation, so to speak) plays an essential role. 
		In the absence of full rotational symmetry, this approach does not readily extend 
		to {Damek–Ricci spaces}—the nonsymmetric generalizations of rank-one 
		Riemannian symmetric spaces of noncompact type. 
		It would be interesting to explore whether analogues of Theorems~\ref{rellich-asymptotic-rss}, ~\ref{rel-main-thm}, and ~\ref{rel-rss-weak-lp} can be established for Damek–Ricci spaces, possibly by developing 
		new harmonic analytic tools that bypass the lack of a classical spherical expansion.
		
		\item 
		In the Baouendi–Grushin setting, Rellich-type inequalities were recently established 
		by Banerjee and Garofalo~\cite{BG2}, building upon ideas of Garofalo–Shen~\cite{GS} and 
		combining Carleman-type estimates, weighted integral inequalities, and the intrinsic 
		homogeneous structure of the operator. 
		A natural direction for further study is to investigate whether analogues of 
		Theorems~\ref{rel-org} and~\ref{BG1-thm} hold in 
		the Heisenberg group, or more generally, in H-type groups.
	\end{enumerate}
	
	\section*{Acknowledgments}The author would like to thank Agnid Banerjee, Rudra Sarkar, and Bernhard Kr\"otz for several helpful 
	discussions, and Swagato K.~Ray for many illuminating conversations and for suggesting 
	interesting directions.  
	The author is supported by the DST--INSPIRE Faculty Fellowship 
	(DST/INSPIRE/04/2024/001099).
	\section*{Data Availability} No data, models, or code were generated or used during the study.
	

\end{document}